\documentclass[twoside,leqno,10pt]{amsart}
\usepackage{amsfonts}
\usepackage{amsmath}
\usepackage{amscd}
\usepackage{amssymb}
\usepackage{amsthm}
\usepackage{amsrefs}
\usepackage{latexsym}
\usepackage{bbm}
\begin{document}

\newtheorem{theorem}[subsection]{Theorem}
\newtheorem{proposition}[subsection]{Proposition}
\newtheorem{lemma}[subsection]{Lemma}
\newtheorem{corollary}[subsection]{Corollary}
\newtheorem{conjecture}[subsection]{Conjecture}
\newtheorem{prop}[subsection]{Proposition}
\numberwithin{equation}{section}
\newcommand{\mr}{\ensuremath{\mathbb R}}
\newcommand{\dif}{\mathrm{d}}
\newcommand{\intz}{\mathbb{Z}}
\newcommand{\ratq}{\mathbb{Q}}
\newcommand{\natn}{\mathbb{N}}
\newcommand{\comc}{\mathbb{C}}
\newcommand{\rear}{\mathbb{R}}
\newcommand{\prip}{\mathbb{P}}
\newcommand{\uph}{\mathbb{H}}
\newcommand{\fief}{\mathbb{F}}
\newcommand{\majorarc}{\mathfrak{M}}
\newcommand{\minorarc}{\mathfrak{m}}
\newcommand{\sings}{\mathfrak{S}}
\newcommand{\fA}{\ensuremath{\mathfrak A}}
\newcommand{\mn}{\ensuremath{\mathbb N}}
\newcommand{\mq}{\ensuremath{\mathbb Q}}
\newcommand{\half}{\tfrac{1}{2}}
\newcommand{\f}{f\times \chi}
\newcommand{\summ}{\mathop{{\sum}^{\star}}}
\newcommand{\chiq}{\chi \bmod q}
\newcommand{\chidb}{\chi \bmod db}
\newcommand{\chid}{\chi \bmod d}
\newcommand{\sym}{\text{sym}^2}
\newcommand{\hhalf}{\tfrac{1}{2}}
\newcommand{\sumstar}{\sideset{}{^*}\sum}
\newcommand{\sumprime}{\sideset{}{'}\sum}
\newcommand{\sumprimeprime}{\sideset{}{''}\sum}
\newcommand{\V}{V\left(\frac{nm}{q^2}\right)}
\newcommand{\sumi}{\mathop{{\sum}^{\dagger}}}
\newcommand{\mz}{\ensuremath{\mathbb Z}}
\newcommand{\leg}[2]{\left(\frac{#1}{#2}\right)}
\newcommand{\muK}{\mu_{\omega}}

\title{Large sieve inequalities for quartic characters}
\date{\today}
\author{Peng Gao and Liangyi Zhao}
\maketitle

\begin{abstract}
In this paper, we prove a large sieve inequality for quartic
Dirichlet characters.  The result is analogous to large sieve
inequalities for the quadratic and cubic Dirichlet characters.
\end{abstract}




\section{Introduction}

The large sieve was an idea originated by J. V. Linnik \cite{JVL1} in 1941 while studying the distribution of quadratic non-residues.  Refinements and extensions in various directions of this idea were made by many \cites{Ba1, Ba2, PrRa, SBLZ, SBLZ3, BD1, Da, DH1, DRHB, Cla, B&Y, PXG, JVL1, Lem, HM2, MVa, Wol, Zha, Zha2, Zha3, Zha4}.  Large sieve results for Dirichlet characters with a fixed order
are particularly useful in analytic number theory. We refer the readers to
\cite{Elliott}, Section 7, for some early large sieve-type results
on general $r$-th order characters.
 Let $(a_n)_{n\in \mathbb{N}}$ be an arbitrary sequence of complex numbers, D. R. Heath-Brown's quadratic large sieve \cite[Theorem 1]{DRHB} states
 that for any $\varepsilon>0$,
\begin{equation}
\label{realfinal}
   \sumstar_{m \leq M} \left| \sumstar_{n \leq N} a_n \Big (\frac {n}{m} \Big ) \right|^2 \ll_{\varepsilon} (MN)^{\varepsilon}(M+N)
   \sumstar_{n \leq N} |a_n|^2,
\end{equation}
   where the asterisks indicate that $m,n$ run over positive odd square-free
   integers and $(\frac {\cdot}{m})$ is the Jacobi symbol. \newline

   Similar to \eqref{realfinal}, Heath-Brown also established the following large sieve inequality involving the cubic symbols \cite[Theorem
   2]{DRHB1}:
\begin{equation}
\label{eq:HBcubic}
 \sumstar_{\substack{m \in \mz[\omega] \\\mathcal{N}(m) \leq M}} \left| \ \sumstar_{\substack{n \in \mz[\omega] \\\mathcal{N}(n) \leq N}} a_n \leg{n}{m}_3 \right|^2
 \ll \left( M + N + (MN)^{2/3} \right)(MN)^{\varepsilon} \sum_{\mathcal{N}(n) \leq N} |a_n|^2,
\end{equation}
where the asterisks indicate that $m,n$ run over square-free elements
of $\mz[\omega], \omega=\exp(2 \pi i/3)$ that are congruent to $1$ modulo 3 and $(\frac {\cdot}{m})_3$ is the cubic residue symbol.  Moreover, here and after, we use $\mathcal{N}(m)$ to denote the norm of $m$. \newline

   Using \eqref{eq:HBcubic}, S. Baier and M. P. Young \cite[Theorem 1.4]{B&Y} proved the following large sieve inequality for cubic Dirichlet characters:
\begin{equation*}
\begin{split}
& \sum\limits_{\substack{Q<q\le 2Q}} \
\sideset{}{^\star}\sum\limits_{\substack{\chi \bmod q\\ \chi^3=\chi_0}}
\left| \ \sumstar\limits_{\substack{M<m\le 2M}} a_m
\chi(m)\right|^2\\
&\ll
(QM)^{\varepsilon}\min\left\{Q^{5/3}+M,Q^{4/3}+Q^{1/2}M,Q^{11/9}+Q^{2/3}M,
Q+Q^{1/3}M^{5/3}+M^{12/5}\right\}
\sumstar\limits_{\substack{M<m\le 2M}}\left| a_m \right|^2,
\end{split}
\end{equation*}
where the star on the sum over $\chi$ restricts the sum to
primitive characters and the asterisks attached to the sum over $m$
indicates that $m$ runs over square-free integers. \newline

It is our goal in this paper to prove a large sieve inequality for quartic
Dirichlet characters.  First we prove the following theorem involving the quartic symbols.
\begin{theorem} \label{mainthm}
Let $M,N$ be positive integers, and let $(a_n)_{n\in \mathbb{N}}$ be an arbitrary
sequence of complex numbers, where $n$ runs over $\mz[i]$. Then we
have
\begin{equation} \label{eq:quartic}
 \sumstar_{\substack{m \in \mz[i] \\\mathcal{N}(m) \leq M}} \left| \ \sumstar_{\substack{n \in \mz[i] \\\mathcal{N}(n) \leq N}} a_n \leg{n}{m}_4 \right|^2
 \ll_{\varepsilon} \left( M + N + (MN)^{3/4} \right)(MN)^{\varepsilon} \sum_{\mathcal{N}(n) \leq N} |a_n|^2,
\end{equation}
   for any $\varepsilon > 0$,where the asterisks indicate that $m$ and $n$ run over square-free elements of $\mz[i]$ that are congruent to $1$ modulo $(1+i)^3$ and $(\frac
{\cdot}{m})_4$ is the quartic residue symbol.
\end{theorem}

   Next we shall establish the following large sieve inequality for quartic Dirichlet characters.
\begin{theorem}
\label{quarticlargesieve} Let $(a_m)_{m\in \mathbb{N}}$ be an
arbitrary sequence of complex numbers. Then
\begin{equation} \label{final}
\begin{split}
&  \sum\limits_{\substack{Q<q\le 2Q}}  \
\sideset{}{^\star}\sum\limits_{\substack{\chi \bmod q\\ \chi^4=\chi_0,
\chi^2 \neq \chi_0}} \left| \ \sumstar\limits_{\substack{M<m\le 2M}}
a_m \chi(m)\right|^2\\ &\ll
(QM)^{\varepsilon}\min\left\{Q^{7/4}+M,Q^{11/8}+Q^{1/2}M, Q^{5/4}+Q^{2/3}M,
Q+Q^{1/2}M+M^{17/7}\right\} \sumstar\limits_{\substack{M<m\le 2M}}\left| a_m
\right|^2,
\end{split}
\end{equation}
where the star on the sum over $\chi$ restricts the sum to
primitive characters and the asterisks attached to the sum over $m$
indicates that $m$ runs over square-free integers.
\end{theorem}

Following the techniques of \cites{DRHB, DRHB1}, Theorem~\ref{mainthm} is proved via recursive uses of the Poisson summation formula.   Theorem~\ref{quarticlargesieve} follows, after some transformations, from Theorem~\ref{mainthm}.   We note that \eqref{eq:quartic} is used in \eqref{square}.   Mark that the characters involved in the second line of \eqref{square} are actually quadratic, since they are squares of the quartic symbol.  Therefore, it is conceivable that the bounds in \eqref{final} can be improved if a large sieve inequality for quadratic characters in $\intz[i]$ is available. \newline

Finally, we wish to mention that it is highly conceivable that these theorems will find applications in the study of families of $L$-functions involving quartic characters, analogous to those results in \cite{B&Y} and \cite{Luo}.

\subsection{Notations} The following notations and conventions are used throughout the paper.\\
\noindent $e(z) = \exp (2 \pi i z) = e^{2 \pi i z}$. \newline
$\widetilde{e} (z) = \exp \left( 2 \pi i ( z + \overline{z} ) \right)$. \newline
$f = O(g)$ or $f \ll g$ means $|f| \leq cg$ for some unspecified
positive constant $c$. \newline

\section{Preliminaries}
\label{sec 2}
\subsection{Quartic symbol and the quartic Gauss sum}
\label{sec2.4}
   The symbol $(\frac{\cdot}{n})_4$ is the quartic
residue symbol in the ring $\mz[i]$.  For a prime $\pi \in \mz[i]$
with $\mathcal{N}(\pi) \neq 2$, the quartic character is defined for $a \in
\mz[i]$, $(a, \pi)=1$ by $\leg{a}{\pi}_4 \equiv
a^{(\mathcal{N}(\pi)-1)/4} \bmod{\pi}$, with $\leg{a}{\pi}_4 \in \{
\pm 1, \pm i \}$. When $\pi | a$, it is defined that
$\leg{a}{\pi}_4 =0$.  Then the quartic character can be extended
to composite $n$ with $(\mathcal{N}(n), 2)=1$ multiplicatively. \newline

 Note that in $\intz[i]$, the ring of cosets modulo $(1+i)^3$ can be represented by $\{0, \pm 1, \pm i, 1+i, 2, 2i \}$ and every ideal coprime to $2$ has a unique
generator congruent to 1 modulo $(1+i)^3$ (\cite[Lemma 7, page
121]{I&R}). Such a generator is called primary. Recall that the
quartic reciprocity law states that for two primary primes  $m, n
\in \mz[i]$,
\begin{equation*}
 \leg{m}{n}_4 = \leg{n}{m}_4(-1)^{((\mathcal{N}(n)-1)/4)((\mathcal{N}(m)-1)/4)}.
\end{equation*}
   Observe that a non-unit
$n=a+bi$ in $\mz[i]$ with $a, b \in \mz$ is congruent to 1 modulo $(1+i)^3$ if and only if $a \equiv 1
\bmod{4}, b \equiv 0 \bmod{4}$ or $a \equiv 3 \bmod{4}, b \equiv 2
\bmod{4}$ by Lemma 6 on page 121 of \cite{I&R}. \newline

Below, we briefly discuss some properties of Gauss sums.  These are dealt with in grater generality in Section 1 of \cite{Patt}.  For $n , r \in \mz[i]$, $n \equiv 1 \bmod {(1+i)^3}$, we set
\begin{equation*}
 g(r,n) = \sum_{x \bmod{n}} \leg{x}{n}_4 \widetilde{e}\left( \frac{rx}{n} \right),
\end{equation*}
where here and after
\begin{equation} \label{etildedef}
\widetilde{e} (z) = \exp \left( 2 \pi i ( z + \overline{z} ) \right).
\end{equation}
   The quartic Gauss sum $g(n)$ is then defined to be $g(n)=g(1,n)$. \newline

   For $(s,n)=1$, we have
\begin{equation*}
 g(rs,n) = \overline{\leg{s}{n}_4} g(r,n).
\end{equation*}
    It's easy to see that the above equality in fact holds for any $s$ when
    $\leg{\cdot}{n}_4$ is a primitive character. \newline

   It's well-known that for square-free $n$'s,
\[ |g(n)|=\sqrt{\mathcal{N}(n)}. \]

Suppose $n \equiv \pm 1 \bmod{{(1+i)^3}}$ with no
rational prime divisor, so $(n, \bar{n}) = 1$. Let $\chi_n$ be a
multiplicative character on $\mz[i]/(n)$, we define
\begin{equation}
\label{tau}
 \tau(\chi_n) =  \sum_{1 \leq x \leq \mathcal{N}(n)} \chi_n(x) e \left( \frac{x}{\mathcal{N}(n)} \right).
\end{equation}
   Now we specify $\chi_n$ to be $\left(\frac{\cdot}{n}\right)_4$. On writing $x = y \bar{n} + \bar{y} n$, where $y$ varies over a set
of representatives in $\mz[i] \bmod{n}$, with $\bar{n}$ being the
complex conjugate of $n$, it's easy to see that
\begin{equation*}
 \tau(\chi_n) =  \sum_{y \bmod{n}} \left(\frac{y \bar{n}}{n}\right)_4 e \left( \frac{y}{n} + \frac{\bar{y}}{\bar{n}} \right) =
 \left(\frac{\bar{n}}{n}\right)_4g(n).
\end{equation*}
   It follows that for $(n_1, n_2)=1$,
\begin{align*}
    \tau(\chi_{n_1n_2})=\left(\frac{\mathcal{N}(n_2)}{n_1}\right)_4\left(\frac{\mathcal{N}(n_1)}{n_2}\right)_4\tau(\chi_{n_1})\tau(\chi_{n_2}),
\end{align*}
    and that if $n$ is square-free
\begin{align*}
  |\tau(\chi_n)|=\sqrt{\mathcal{N}(n)}.
\end{align*}
  Similarly, we have for $n$ square-free
\begin{align*}
  |\tau(\chi^2_n)|=\sqrt{\mathcal{N}(n)}.
\end{align*}

\subsection{Primitive quartic Dirichlet characters}
\label{sec2.5}
   The classification of all the primitive cubic characters of conductor $q$ coprime to $3$ is given in
   \cite{B&Y}.
Similarly, one can give a classification of all the primitive
quartic characters of conductor $q$ coprime to
    $2$.  Every such character is of the form $m \rightarrow
(\frac{m}{n})_4$ for some $n \in \mz[i]$, with $n \equiv 1
\bmod{(1+i)^3}$, $n$ square-free and not divisible by any rational
primes and $\mathcal{N}(n) = q$.

\section{Strategy for the proof of Theorem \ref{mainthm}}

   Our proof of Theorem \ref{mainthm} uses the ideas in \cites{DRHB, DRHB1}. We first
   estimate
\begin{equation*}
   {\sum}_1=\sumstar_{\substack{m \in \mz[i] \\M < \mathcal{N}(m) \leq 2M}} \left| \ \sumstar_{\substack{n \in \mz[i] \\N < \mathcal{N}(n) \leq 2N}} a_n \leg{n}{m}_4
   \right|^2.
\end{equation*}
   We further simplify notation by supposing that
the coefficients $a_n$ are supported on such integers $n \in
\mz[i]$ satisfying $N < \mathcal{N}(n)  \leq 2N$.
   We begin by defining the norm
\begin{equation*}
   \mathcal{B}_1(M,N)=\sup \left\{ {\sum}_1: \sum_n|a_n|^2=1 \right\}.
\end{equation*}
Therefore, we need to show
\begin{equation*}
   \mathcal{B}_1(M,N) \ll_{\varepsilon} (MN)^{\varepsilon} \left( M + N + (MN)^{3/4} \right).
\end{equation*}
Introducing a smooth weight function, we have
\begin{equation*}
   {\sum}_1 \ll \sum_m \exp \left (-2\pi \frac{\mathcal{N}(m)}{M} \right) \left| \ \sumstar_{\substack{n \in \mz[i] \\ N <\mathcal{N}(n) \leq 2N}} a_n \leg{n}{m}_4
   \right|^2,
\end{equation*}
   the sum being over all $m \in \mz[i]$ for which $m \equiv 1 \bmod {(1+i)^3}$. If we now expand
the above expression we obtain sums of the form
\begin{equation} \label{2.1}
   \sum_m \exp \left( -2\pi \frac{\mathcal{N}(m)}{M} \right) \leg{n_1}{m}_4\overline{\leg{n_2}{m}}_4.
\end{equation}
   We note the following analogue of Lemma 2 of Heath-Brown and Patterson
   \cite{H&P}.  As the proof is similar, we omit it here.
\begin{lemma}
\label{lem1} Let $\chi$ be a character of modulus $f \neq 1$, not
necessarily primitive. Then, for $w \leq 1, \varepsilon>0$,
\begin{equation*}
  \theta(w, \chi)=\sum_{\substack{ a \equiv 1 \bmod{ (1+i)^3} \\ (a,f)=1}} \chi(a)e^{-2\pi \mathcal{N}(a)w} \ll E(\chi)w^{-1}+\mathcal{N}(f)^{1/2+\varepsilon},
\end{equation*}
  where $E(\chi)=1$ if $\chi$ is principal, $0$ otherwise. The
  implied constant depends only on $\varepsilon$.
\end{lemma}

Lemma \ref{lem1} implies that each of these sums in \eqref{2.1}
are $O\left( \mathcal{N}(n_1n_2)^{1/2+\varepsilon} \right)$, provided that the
character involved is non-principal. Since $n_1$ and $n_2$ are
square-free, $\leg{n_1}{m}_4\overline{\leg{n_2}{m}_4}$ is principle only if  $n_1 = n_2$.
It follows that
\[ {\sum}_1 \ll_{\varepsilon}
   N^{\varepsilon} \left( M\sum_n|a_n|^2+N\sum_{n_1,n_2}|a_{n_1}a_{n_2}| \right)  \ll_{\varepsilon}
   N^{\varepsilon} \left( M+N^2 \right)\sum_n|a_n|^2. \]
 We therefore have
\begin{equation}
\label{initialest}
   \mathcal{B}_1(M,N) \ll_{\varepsilon} N^{\varepsilon}\left( M+N^2 \right).
\end{equation}
    This will be the starting point for an iterative bound for
    $\mathcal{B}_1(M,N)$. \newline

   Similar to the proof of \cite[Lemma 1]{DRHB}, using the duality
   principle (see for example, \cite[Chap. 9]{HM}) and the quartic
   reciprocity law by considering the case for $n=a+bi$ with
   $a \equiv 1 \bmod{4}, b \equiv 0 \bmod{4}$ or $a
\equiv 3 \bmod{4}, b \equiv 2 \bmod{4}$ (and similarly for $m$),
we can establish the following lemma.
\begin{lemma}
\label{lem2} We have $\mathcal{B}_1(M,N) \leq 2
\mathcal{B}_1(N,M)$. Moreover, there exist coefficients $a'_n,
a''_n$ with $|a'_n|=|a''_n|=|a_n|$ such that
\begin{equation*}
\begin{split}
\sumstar_{\substack{m \in \mz[i] \\M <\mathcal{N}(m) \leq 2M}} \left| \
\sumstar_{\substack{n \in \mz[i] \\N < \mathcal{N}(n) \leq 2N}} a_n
\leg{n}{m}_4
   \right|^2  & \leq 2 \sumstar_{\substack{m \in \mz[i] \\M< \mathcal{N}(m) \leq 2M}} \left| \ \sumstar_{\substack{n \in \mz[i] \\N < \mathcal{N}(n) \leq 2N}} a'_n \leg{m}{n}_4
   \right|^2 \\
   & \leq 4 \sumstar_{\substack{m \in \mz[i] \\M< \mathcal{N}(m) \leq 2M}} \left| \ \sumstar_{\substack{n \in \mz[i] \\N < \mathcal{N}(n) \leq 2N}} a''_n \leg{n}{m}_4
   \right|^2.
\end{split}
\end{equation*}
\end{lemma}

    Our next lemma is a trivial modification of Lemma 9 of \cite{DRHB}, which shows that the norm $\mathcal{B}_1(M,N)$ is essentially
    increasing.
\begin{lemma}
\label{lem3} There is an absolute constant $C > 0$ as follows.
Let $M_1, N \geq 1$ and $M_2 \geq CM_1\log(2M_1N)$. Then
\begin{equation*}
   \mathcal{B}_1(M_1,N)\leq C \mathcal{B}_1(M_2,N).
\end{equation*}
   Similarly, if $M, N_1 \geq 1$ and $N_2 \geq CN_1\log(2N_1M)$. Then
\begin{equation*}
   \mathcal{B}_1(M,N_1)\leq C \mathcal{B}_1(M,N_2).
\end{equation*}
\end{lemma}

Next, we define
\begin{equation*}
   \mathcal{B}_2(M,N)=\sup \left\{ {\sum}_2: \sum_n|a_n|^2=1 \right\},
\end{equation*}
    where
\begin{equation}
\label{sum2}
   {\sum}_2=\sum_{\substack{m \in \mz[i] \\M < \mathcal{N}(m) \leq 2M}} \left| \ \sumstar_{\substack{n \in \mz[i] \\N < \mathcal{N}(n) \leq 2N}} a_n \leg{m}{n}_4
   \right|^2,
\end{equation}
   the summation over $m$ running over all integers of $\mz[i]$ in the relevant range. \newline

   It follows directly from Lemma \ref{lem2} that
\begin{equation}
\label{12comparison}
  \mathcal{B}_1(M,N) \leq 2\mathcal{B}_2(M,N).
\end{equation}
For the other direction, we have the following.
\begin{lemma}
\label{lem4} There exist $X, Y \gg 1$ such that $XY^3 \ll M$ and
\begin{equation*}
   \mathcal{B}_2(M,N) \ll (\log M)^3M^{1/2}X^{-1/2}Y^{-3/2}\min(Y\mathcal{B}_1(X,N), X\mathcal{B}_1(Y,N)).
\end{equation*}
\end{lemma}
\begin{proof}
   To handle $\sum_2$ we write each of the integers $m$ occurring in the outer summation
of \eqref{sum2} in the form $m = ab^2c^3d$, where $a, b, c\equiv 1
\bmod {(1+i)^3}$ are square-free, and $d$ is a product of a unit,
a power of $1+i$, and a fourth power (so that $d$ can be written
as $d=u(1+i)^je^4$ where $u$ is a unit, $0 \leq j \leq 3$ and $e
\in \mz[i]$). We split the available ranges for $a, b, c$ and $d$
into sets $X < \mathcal{N}(a) \leq 2X, Y < \mathcal{N}(b) \leq 2Y, Z < \mathcal{N}(c) \leq 2Z$
and $W < \mathcal{N}(d) \leq 2W$, where $X, Y, Z$ and $W$ are powers of $2$.
There will therefore be $O(\log^3M)$ possible quadruples $X, Y,Z,
W$. We may now write
\begin{equation*}
   {\sum}_2 \ll \sum_{X,Y,Z, W}{\sum}_2(X,Y,Z,W)
\end{equation*}
  accordingly, so that
\begin{equation*}
   {\sum}_2 \ll (\log^3M){\sum}_2(X,Y,Z,W)
\end{equation*}
   for some quadruple $X,Y,Z, W$. However,
\begin{equation*}
   {\sum}_2(X,Y,Z,W) \leq \sum_{b,c,d}\sumstar_{\substack{a  \in \mz[i] \\X'<\mathcal{N}(a) \leq 2X'}}
   \left| \ \sumstar_{\substack{n \in \mz[i] \\N < \mathcal{N}(n) \leq 2N}} a_n \leg{b^2c^3d}{n}_4 \leg{a}{n}_4
   \right|^2,
\end{equation*}
   where $X'=X'(b,c, d)=M/\mathcal{N}(b^2c^3d)$. It is easy to see that $X/2 \leq X' \leq
   2X$, and hence by Lemma \ref{lem2}
\[ {\sum}_2(X,Y,Z, W) \ll \sum_{b,c,d}\mathcal{B}_1(X',N)\sum_{n}
   |a_n|^2 \ll YZW^{1/4}\max \left\{ \mathcal{B}_1(X',N): X/2 \leq X' \leq 2X \right\} \sum_{n}
   |a_n|^2, \]
 since there are $O(W^{1/4})$ possible integers $d$. \newline

   In the same way we have
\begin{align*}
   {\sum}_2(X,Y,Z, W) &\leq \sum_{a, b, d}\sumstar_{\substack{c  \in \mz[i] \\Z'<\mathcal{N}(c) \leq 2Z'}}
   \left| \ \sumstar_{\substack{n \in \mz[i] \\N < \mathcal{N}(n) \leq 2N}} a_n \leg {ab^2d}{n}_4 \leg {c^3}{n}_4
   \right|^2 \\
   &= \sum_{a,b, d}\sumstar_{\substack{c  \in \mz[i] \\Z'<\mathcal{N}(c) \leq 2Z'}}
   \left| \ \sumstar_{\substack{n \in \mz[i] \\N < \mathcal{N}(n) \leq 2N}} \overline{a}_n \overline{\leg {ab^2d}{n}}_4 \overline{\leg {c^3}{n}}_4
   \right|^2 \\
   &= \sum_{a,b,d}\sumstar_{\substack{c  \in \mz[i] \\Z'<\mathcal{N}(c) \leq 2Z'}}
   \left| \ \sumstar_{\substack{n \in \mz[i] \\N < \mathcal{N}(n) \leq 2N}} \overline{a}_n \overline{\leg {ab^2d}{n}}_4 \leg {c}{n}_4
   \right|^2 \\
   & \ll \ \sum_{a,b,d}\mathcal{B}_1(Z',N)\sum_{n}
   |a_n|^2 \\
   & \ll XYW^{1/4}\max \left\{ \mathcal{B}_1(Z',N): Z \ll Z' \ll Z \right\} \sum_{n}
   |a_n|^2,
\end{align*}
   where $Z'=Z'(a,b,d)=M/\mathcal{N}(ab^2d)$. As $Y \ll M^{1/2}X^{-1/2}Z^{-3/2}W^{-1/2}$, we see that
\begin{equation*}
   \mathcal{B}_2(M,N) \ll (\log M)^3M^{1/2}X^{-1/2}Z^{-3/2}W^{-1/4}\min(Z\mathcal{B}_1(X,N), X\mathcal{B}_1(Z,N)).
\end{equation*}
   The assertion of the lemma now follows on replacing $Z$ by $Y$ above.
\end{proof}

   As in \cite{DRHB}, we introduce an infinitely differentiable weight function $W: \mr \rightarrow \mr$, defined by
\begin{equation}
\label{W}
W(x)= \begin{cases} \exp\left(\frac{-1}{(2x-1)(5-2x)}\right),
\qquad & \text{if }  \frac{1}{2}<x<\frac{5}{2},\\  0, \qquad  &
\text{otherwise}. \end{cases}
\end{equation}
We now have
\begin{equation}
\label{3.7}
   {\sum}_2 \ll \sum_{m \in \mz[i]}W\left( \frac {\mathcal{N}(m)}{M}\right) \left| \sum_{\substack{n \in \mz[i] \\N < \mathcal{N}(n) \leq 2N}} a_n \leg {m}{n}_4
   \right|^2,
\end{equation}
   where we recall that we can drop the conditions on $a_n$ on the inner sum above of the right-hand side expression by supposing that
the coefficients $a_n$ are supported on square-free integers $n
\equiv 1 \bmod {(1+i)^3} \in \mz[i]$ lying in the range $N < \mathcal{N}(n)
\leq 2N$. \newline

Expanding the sum on the right-hand side of \eqref{3.7}, we obtain
\begin{equation*}
   {\sum}_2 \ll \sum_{n_1,n_2}a_{n_1}\overline{a}_{n_2}\sum_{m \in \mz[i]}
   W\left(\frac {\mathcal{N}(m)}{M}\right)
   \leg{m}{n_1}_4\overline{\leg{m}{n_2}}_4.
\end{equation*}
   We set
\begin{equation*}
   {\sum}_3={\sum}_3(M,N)= \sum_{(n_1,n_2)=1}a_{n_1}\overline{a}_{n_2}\sum_{m \in \mz[i]}
   W\left(\frac {\mathcal{N}(m)}{M}\right)
   \leg{m}{n_1}_4\overline{\leg{m}{n_2}}_4
\end{equation*}
   and define
\begin{equation} \label{B3def}
   \mathcal{B}_3(M,N)=\sup \left\{ {\sum}_3: \sum_n|a_n|^2=1 \right\}.
\end{equation}
    Similar to \cite[Lemma 7]{DRHB, DRHB1}, we have the
    following
\begin{lemma}
\label{lem5} Let $\varepsilon > 0$ be given. Then there exist
positive integers $\Delta_2 \geq \Delta_1$ such that
\begin{equation*}
   \mathcal{B}_2(M,N) \ll_{\varepsilon} N^{\varepsilon} \mathcal{B}_3 \left(\frac {M}{\Delta_1}, \frac {N}{\Delta_2}\right).
\end{equation*}
\end{lemma}

   We complete the chain of relations amongst the various norms by giving the following estimate for $\mathcal{B}_3(M,N)$ in terms of
   $\mathcal{B}_2(M,N)$.
\begin{lemma}
\label{lem6} Let $N \geq 1$. Then for any $\varepsilon>0$ we have
\begin{equation*}
   \mathcal{B}_3(M,N) \ll_{\varepsilon} MN^{4\varepsilon-1}\max \left\{ \mathcal{B}_2(K,
   N): K \leq N^2/M
   \right\}+M^{-1}N^{3+4\varepsilon}\sum_{K>N^2/M}K^{-2-\varepsilon}\mathcal{B}_2(K,
   N),
\end{equation*}
   where $K$ runs over powers of $2$.
\end{lemma}
This bound uses the Poisson summation formula and is the key
in the proof of Theorem \ref{mainthm}.  Note that it does not cover the case in which $N =1/2$, say, for which we have the trivial bound
\begin{equation}
\label{trivialbound}
   \mathcal{B}_3(M,N) \ll_{\varepsilon} M, \hspace{0.1in} (N \leq 1).
\end{equation}
Section \ref{sec4} will be devoted to the proof of Lemma~\ref{lem6}.

\section{Proof of Lemma \ref{lem6}} \label{sec4}

Our proof of Lemma \ref{lem6} requires the application of the Poisson
summation formula. We shall write
\begin{equation*}
   \chi(m)=\leg{m}{n_1}_4\overline{\leg{m}{n_2}}_4,
\end{equation*}
which is a primitive character (on the group
$(\mz[i]/(n_1n_2))^{\times}$) to modulus $q = n_1n_2$, provided
that $n_1$, $n_2$ and 2 are pair-wise coprime and that $n_1$ and $n_2$ are
square-free.
\begin{lemma}
\label{lem7} With the above notations we have
\begin{equation*}
   \sum_{m \in \mz[i]}W\left(\frac {\mathcal{N}(m)}{M}\right)\chi(m)=\frac {\chi(-2i)g(n_1)\overline{g(n_2)}M}{\mathcal{N}(q)}\leg{n_2}{n_1}_4\overline{\leg{n_1}{n_2}}_4\leg{-1}{n_2}_4
   \sum_{m \in \mz[i]}\widetilde{W}\left(\sqrt{\frac {\mathcal{N}(k)M}{\mathcal{N}(q)}}\right)\overline{\chi}(k),
\end{equation*}
   where
\begin{equation*}
   \widetilde{W}(t)=\int\limits^{\infty}_{-\infty}\int\limits^{\infty}_{-\infty}W(\mathcal{N}(x+yi))\widetilde{e}\left( \frac{t(x+yi)}{2i}\right)\dif x  \dif y,
\end{equation*}
   for non-negative $t$. Here $\widetilde{e}(z)$ is defined in \eqref{etildedef} and $g(n)$ is the Gauss sums defined in Section \ref{sec2.4}.
\end{lemma}
\begin{proof}
This lemma is analogous to Lemma 10 in \cite{DRHB1} and the proof is very similar.  The differences include we need to start with the Poisson summation formula for
   $\mz[i]$, which takes the form.
\begin{equation*}
   \sum_{j \in \mz[i]}f(j)=\sum_{k \in
   \mz[i]}\int\limits^{\infty}_{-\infty}\int\limits^{\infty}_{-\infty}f(x+yi)\widetilde{e}\left( \frac {k(x+yi)}{2i} \right)\dif x  \dif y.
\end{equation*}
We omit the details of the rest of proof as it simply goes along the same lines as the proof of Lemma 10 in \cite{DRHB1}.
\end{proof}

    Our next result will be used to separate the variables in a function of a
product, which is Lemma 12 of \cite{DRHB}.
\begin{lemma}
\label{lem8} Let $\rho: \mr \rightarrow \mr$ be an infinitely
differentiable function whose derivatives satisfy the bound
\[ \rho{(k)}(x) \ll_{k,A}|x|^{-A} \]
for $|x| \geq 1$, for any positive constant
$A$. Let
\begin{equation*}
   \rho_{+}(s)=\int\limits^{\infty}_{0}\rho(x)x^{s-1} \dif x,
   \hspace{0.1in} \rho_{-}(s)=\int\limits^{\infty}_{0}\rho(-x)x^{s-1} \dif
   x.
\end{equation*}
   Then $\rho_{+}(s)$ and $\rho_{-}(s)$ are holomorphic in $\Re(s) = \sigma > 0$, and
   satisfy
\begin{equation*}
   \rho_{+}(s), \ \rho_{-}(s) \ll_{A, \sigma} |s|^{-A},
\end{equation*}
in that same domain, for any positive constant $A$. Moreover if $\sigma> 0$ we have
\begin{equation*}
   \rho(x)=\frac {1}{2\pi
   i}\int\limits^{\sigma+i\infty}_{\sigma-i\infty}\rho_{+}(s)x^{-s} \dif s \; \; \; \mbox{and} \; \; \;
   \rho(-x)=\frac {1}{2\pi
   i}\int\limits^{\sigma+i\infty}_{\sigma-i\infty}\rho_{-}(s)x^{-s} \dif s
\end{equation*}
   for any positive $x$.
\end{lemma}

We are now ready to prove Lemma~\ref{lem6}.

\begin{proof} [Proof of Lemma~\ref{lem6}]
   In the notation of Lemma \ref{lem7} we have
\begin{equation*}
   {\sum}_3(M,N)=\sum_{(n_1,n_2)=1}a_{n_1}\overline{a}_{n_2}\sum_{m \in \mz[i]}W\left(\frac {\mathcal{N}(m)}{M}\right)\chi(m).
\end{equation*}
   We proceed to evaluate the inner sum using Lemma \ref{lem7}, whence
\begin{equation}
\label{5.1}
   {\sum}_3(M,N)=M\sum_{k \in \mz[i]}\sum_{(n_1,n_2)=1}c_{n_1}\overline{c}_{n_2}\leg{n_2}{n_1}_4\overline{\leg{n_1}{n_2}}_4\leg{-1}{n_2}_4\widetilde{W}\left(\sqrt{\frac {\mathcal{N}(k)M}{\mathcal{N}(n_1n_2)}}\right)\overline{\chi}(k),
\end{equation}
   where
\begin{equation*}
    c_n=a_n\leg{-2i}{n}_4\frac {g(n)}{\mathcal{N}(n)}.
\end{equation*}
   Note by the law of quartic reciprocity, we have
\begin{equation*}
   \leg{n_2}{n_1}_4\overline{\leg{n_1}{n_2}}_4=(-1)^{((\mathcal{N}(n_1)-1)/4)((\mathcal{N}(n_2)-1)/4)}.
\end{equation*}
    Now we let
\[ S_1 =\{n \in \mz[i]: N<\mathcal{N}(n)\leq 2N, n \hspace{0.05in} \text{square-free}, n=a+bi, a, b \in \mz, a \equiv 1 \bmod{4}, b \equiv 0 \bmod{4} \},  \]
and
\[ S_2 =\{n \in \mz[i]: N<\mathcal{N}(n)\leq 2N, n \hspace{0.05in} \text{square-free}, n=a+bi, a, b \in \mz, a \equiv 3 \bmod{4}, b \equiv 2 \bmod{4} \}. \]
   We can then recast the inner sum in \eqref{5.1} as
\begin{align*}
   &\sum_{(n_1,n_2)=1} \cdots \\
   =& \sum_{\substack{(n_1,n_2)=1 \\ n_1 \in S_1, n_2 \in S_1}} \cdots
   \; \; +\sum_{\substack{(n_1,n_2)=1 \\ n_1 \in S_1, n_2 \in S_2}} \cdots \; \; +\sum_{\substack{(n_1,n_2)=1 \\ n_1 \in S_2, n_2 \in S_1}} \cdots \; \;
   +\sum_{\substack{(n_1,n_2)=1 \\ n_1 \in S_2, n_2 \in S_2}} \cdots \; \; -2\sum_{\substack{(n_1,n_2)=1 \\ n_1 \in S_2, n_2 \in S_2}} \cdots \\
   =& \sum_{(n_1,n_2)=1}c_{n_1}\overline{c}_{n_2}\leg{-1}{n_2}_4\widetilde{W}\left(\sqrt{\frac {\mathcal{N}(k)M}{\mathcal{N}(n_1n_2)}}\right)\overline{\chi}(k)-2\sum_{(n_1,n_2)=1}c'_{n_1}\overline{c'}_{n_2}\leg{-1}{n_2}_4\widetilde{W}\left(\sqrt{\frac {\mathcal{N}(k)M}{\mathcal{N}(n_1n_2)}}\right)\overline{\chi}(k),
\end{align*}
   where we let $c'_n=c_n$ if $n \in S_2$ and $0$ otherwise. Due
   to similarities, it suffices to estimate
\begin{equation*}
   M\sum_{k \in
\mz[i]}\sum_{(n_1,n_2)=1}c_{n_1}\overline{c}_{n_2}\leg{-1}{n_2}_4\widetilde{W}\left(\sqrt{\frac
{\mathcal{N}(k)M}{\mathcal{N}(n_1n_2)}}\right)\overline{\chi}(k),
\end{equation*}
   Note that $k = 0$ may be omitted if $N \geq 1$, since then $\mathcal{N}(n_1n_2) > 1$ and $\chi(0) = 0$,
the character being non-trivial. We may now apply Lemma \ref{lem8}
to the function $\rho(x) = \widetilde{W}(x)$, which satisfies the necessary conditions of the lemma, as one sees by repeated
integration by parts. We decompose the available $k$ into sets for
which $K < \mathcal{N}(k) \leq 2K$, where $K$ runs over powers of $2$, and
use
\begin{equation} \label{sigmadef}
\sigma = \left\{ \begin{array}{ll} \varepsilon , &  \mbox{for} \;  K \leq N^2/M, \\  4 +
\varepsilon , &  \mbox{otherwise}.
\end{array} \right.
\end{equation}
This gives
\begin{align*}
    {\sum}_3 & \ll_{\varepsilon} M
    \sum_{K}(KM)^{-\sigma/2}\int\limits^{\infty}_{-\infty}|\rho_{+}(\sigma+it)||S(\sigma+it)| \dif t,
\end{align*}
  where
\[
   S(s)=\sum_{K < \mathcal{N}(k) \leq 2K}\left| \sum_{(n_1,n_2)=1}b_{n_1}b'_{n_2}\leg{-1}{n_2}_4\overline{\chi}(k)\right|, \; \; \;  \mbox{with} \; \; \;
    b_n=c_n\mathcal{N}(n)^{s/2} \; \; \mbox{and} \; \; b'_n=\overline{c}_n\mathcal{N}(n)^{s/2} .
\]
We use the M\"obius function to detect the coprimality condition in the inner sum of $S(s)$, giving
\begin{align*}
 S(s)  &\ll \sum_{d}\sum_{K < \mathcal{N}(k) \leq 2K}\left| \sum_{d|(n_1,n_2)}b_{n_1}b'_{n_2}\leg{-1}{n_2}_4\overline{\chi}(k)\right| \\
&= \sum_{d}\sum_{K < \mathcal{N}(k) \leq 2K}\left|
\sum_{d|n}b_{n}\overline{\leg{k}{n}}_4\right| \left|
\sum_{d|n}b'_{n}\leg{-k}{n}_4\right| \leq S^{1/2}_1S^{1/2}_2,
\end{align*}
   by Cauchy's inequality, where
\[  S_1 = \sum_{d}\sum_{K < \mathcal{N}(k) \leq 2K}\left| \sum_{d|n}b_{n}\overline{\leg{k}{n}}_4\right|^2 \]
and satisfies the bound
\[  S_1 \leq \sum_{d}\mathcal{B}_2(K,N)\sum_{d|n}|b_n|^2  \leq \mathcal{B}_2(K,N)\sum_{n}d(n)|a_n|^2\mathcal{N}(n)^{\sigma-1} \ll_{\varepsilon} N^{\varepsilon+\sigma-1}\mathcal{B}_2(K,N). \]
$S_2$ can be treated similarly. It follows then that
\begin{equation*}
   S(s) \ll_{\varepsilon} N^{\varepsilon+\sigma-1}\mathcal{B}_2(K,N),
\end{equation*}
   and since
\begin{align*}
     \int\limits^{\infty}_{-\infty}|\rho_{+}(\sigma+it)| \dif t \ll_{\varepsilon} 1,
\end{align*}
we infer, mindful of our choices of $\sigma$ in \eqref{sigmadef}, that
\begin{align*}
    {\sum}_3 & \ll_{\varepsilon}  MN^{4\varepsilon-1}\max \left\{ \mathcal{B}_2(K,
   N): K \leq N^2/M
   \right\}+M^{-1}N^{3+4\varepsilon}\sum_{K>N^2/M}K^{-2-\varepsilon}\mathcal{B}_2(K,
   N),
\end{align*}
Recalling the definition of $\mathcal{B}_3(M,N)$ in \eqref{B3def}, we have completed the proof of Lemma \ref{lem6}.
   \end{proof}

\section{The Recursive Estimate and the Proof of Theorem~\ref{mainthm}}

    Lemmas \ref{lem4}, \ref{lem5} and \ref{lem6} allow us to estimate $\mathcal{B}_1(M,N)$ recursively, as follows.
\begin{lemma}
\label{lem9} Suppose that $3/2 < \xi \leq 2$, and that
\begin{equation}
\label{6.1}
   \mathcal{B}_1(M,N)\ll_{\varepsilon} (MN)^{\varepsilon} \left( M+N^{\xi} + (MN)^{3/4} \right).
\end{equation}
   for any $\varepsilon > 0$. Then
\begin{equation*}
   \mathcal{B}_1(M,N)\ll_{\varepsilon} (MN)^{\varepsilon} \left( M+N^{(9\xi-6)/(4\xi-1)} + (MN)^{3/4} \right).
\end{equation*}
   for any $\varepsilon > 0$.
\end{lemma}
\begin{proof}
  By the symmetry expressed in Lemma \ref{lem2} the hypothesis
\eqref{6.1} yields
\begin{align*}
     \mathcal{B}_1(M,N)\ll_{\varepsilon} \left( M^{\xi} + N + (MN)^{3/4}\right) (MN)^{\varepsilon}.
\end{align*}
    It follows from \eqref{initialest} that the above estimation is valid with
    $\xi=2$. We now feed this into Lemma \ref{lem4}, whence
\begin{equation} \label{eqwithmin}
   \mathcal{B}_2(M,N) \ll_{\varepsilon} (MN)^{2\varepsilon}M^{1/2}X^{-1/2}Y^{-3/2}\min(Yf(X,N),
   Xf(Y,N)),
\end{equation}
   where
\begin{equation*}
     f(Z,N)=Z^{\xi} + N + (ZN)^{3/4}.
\end{equation*}
   If $X \geq Y$ we bound the minimum in \eqref{eqwithmin} by $Yf(X,N)$, whence
\begin{equation*}
   \mathcal{B}_2(M,N) \ll_{\varepsilon} (MN)^{2\varepsilon}M^{1/2}X^{-1/2}Y^{-3/2}\left( YX^{\xi} + YN +
   Y(XN)^{3/4} \right).
\end{equation*}
   Here we have
\begin{equation*}
   M^{1/2}X^{-1/2}Y^{-3/2}YX^{\xi} \ll M^{\xi}Y^{1-3\xi}
\end{equation*}
   since $X \ll MY^{-3}$. On recalling that $\xi > 3/2 >1/3$ and $Y \gg 1$ we see that
this is $O(M^{\xi})$. Moreover
\begin{equation*}
   M^{1/2}X^{-1/2}Y^{-3/2}YN=M^{1/2}X^{-1/2}Y^{-1/2}N \ll
   M^{1/2}N.
\end{equation*}
   Finally
\[  M^{1/2}X^{-1/2}Y^{-3/2}Y(XN)^{3/4}=M^{1/2}X^{1/4}Y^{-1/2}N^{3/4}
   \ll M^{3/4}N^{3/4}
   \ll M^{1/2}N+M^{3/2} \leq
   M^{1/2}N+M^{\xi}, \]
   since $\xi > 3/2$. It follows that
\begin{equation}
\label{B2bound}
   \mathcal{B}_2(M,N) \ll_{\varepsilon}
   (MN)^{2\varepsilon} \left( M^{1/2}N+M^{\xi} \right)
\end{equation}
   when $X \geq Y$. \newline

   In the alternative case we bound the minimum in \eqref{eqwithmin} by $Xf(Y,N)$, whence
\begin{equation*}
   \mathcal{B}_2(M,N) \ll_{\varepsilon} (MN)^{2\varepsilon}M^{1/2}X^{-1/2}Y^{-3/2} \left( XY^{\xi} + XN + X(YN)^{3/4} \right).
\end{equation*}
   Here
\begin{equation*}
   M^{1/2}X^{-1/2}Y^{-3/2}XY^{\xi} \ll M^{1/2}X^{1/2}Y^{1/2} \ll M
   \ll M^{\xi}
\end{equation*}
   since $\xi \leq 2$ and $XY \ll M$. Moreover
\begin{equation*}
   M^{1/2}X^{-1/2}Y^{-3/2}XN  = M^{1/2}X^{1/2}Y^{-3/2}N \ll
   M^{1/2}N
\end{equation*}
  since we are now supposing that $Y \geq X$.
  Finally
\[ M^{1/2}X^{-1/2}Y^{-3/2}X(YN)^{3/4}=M^{1/2}X^{1/2}Y^{-3/4}N^{3/4} \ll
   M^{1/2}Y^{-1/4}N^{3/4} \ll  M^{1/2}N^{3/4} \ll M^{1/2}N+M^{\xi}, \]
   as before. It follows that \eqref{B2bound} holds when $Y \geq X$ too. It will be convenient to
observe that \eqref{B2bound} still holds when $M < 1/2$ , since
then $\mathcal{B}_2(M,N)= 0$. \newline

   We are now ready to use \eqref{B2bound} (with a new value for $\varepsilon$) in Lemma \ref{lem6}, to obtain
a bound for $\mathcal{B}_3(M,N)$. We readily see that
\begin{equation*}
   \max\left\{ \mathcal{B}_2(K,
   N): K \leq N^2/M
   \right\} \ll_{\varepsilon} N^{\varepsilon}\left( M^{-1/2}N^2+M^{-\xi}N^{2\xi} \right)
\end{equation*}
   and
\begin{equation*}
   \sum_{K>N^2/M}K^{-2-\varepsilon}\mathcal{B}_2(K,
   N) \ll_{\varepsilon}
   N^{\varepsilon}\left( M^{3/2}N^{-2}+M^{2-\xi}N^{2\xi-4} \right).
\end{equation*}
   Thus, if $N \geq 1$, we will have
\begin{equation*}
    \mathcal{B}_3(M,
   N) \ll_{\varepsilon}
   N^{5\varepsilon}\left( M^{1/2}N+M^{1-\xi}N^{2\xi-1} \right).
\end{equation*}
   When this is used in Lemma \ref{lem5} we find that when $N/\Delta_2 \geq
   1$,
\begin{align*}
    \mathcal{B}_3\left( \frac {M}{\Delta_1}, \frac {N}{\Delta_2} \right) &\ll_{\varepsilon}
   N^{5\varepsilon}\left( M^{1/2}N+M^{1-\xi}N^{2\xi-1}\Delta^{\xi-1}_1\Delta^{1-2\xi}_2 \right) \leq N^{5\varepsilon} \left( M^{1/2}N+M^{1-\xi}N^{2\xi-1}\Delta^{-\xi}_2 \right) \\
   & \leq N^{5\varepsilon}\left( M^{1/2}N+M^{1-\xi}N^{2\xi-1} \right).
\end{align*}
   Note that when $M \geq N$, we have $M^{1/2}N \leq (MN)^{3/4}$
   and when $M \leq N$, we have $(MN)^{3/4} \leq
   M^{1-\xi}N^{2\xi-1}$ since $\xi > 3/2$. Thus we conclude that
\begin{equation*}
    \mathcal{B}_2(M,
   N) \ll_{\varepsilon}
   N^{6\varepsilon}\left( (MN)^{3/4}+M^{1-\xi}N^{2\xi-1} \right),
\end{equation*}
   provided that $N/\Delta_2 \geq 1$. In the alternative case \eqref{trivialbound} applies, whence
\begin{equation*}
    \mathcal{B}_2(M,
   N) \ll_{\varepsilon}
   (MN)^{6\varepsilon}\left( M+(MN)^{3/4}+M^{1-\xi}N^{2\xi-1} \right),
\end{equation*}
   In view of Lemma \ref{lem3} and \eqref{12comparison} we may now deduce that
\begin{align*}
    \mathcal{B}_1(M,
   N)  \leq \mathcal{B}_1(M',
   N) \ll \mathcal{B}_2(M',
   N) \ll_{\varepsilon}
   (M'N)^{6\varepsilon}\left( M'+(M'N)^{3/4}+{M'}^{1-\xi}N^{2\xi-1} \right) ,
\end{align*}
  for any $M' \geq CM \log(2MN)$. Note that when ${M}^{4\xi-1}
  \leq N^{8\xi-7}$, we have
\begin{align*}
    (MN)^{3/4} \leq {M}^{1-\xi}N^{2\xi-1}.
\end{align*}
   We shall now choose
\begin{align*}
    M'=C\max \left\{ M, N^{(8\xi-7)/(4\xi-1)} \right\} \log(2MN) ,
\end{align*}
   so that when $M \geq N^{(8\xi-7)/(4\xi-1)}$, we have
\begin{align*}
     M'+(M'N)^{3/4}+{M'}^{1-\xi}N^{2\xi-1} \ll (MN)^{\varepsilon} \left( M+(MN)^{3/4} \right) ,
\end{align*}
  while when $M \leq N^{(8\xi-7)/(4\xi-1)}$, we have
\begin{align*}
    M'+(M'N)^{3/4}+{M'}^{1-\xi}N^{2\xi-1} &
     \ll (MN)^{\varepsilon} \left( N^{(8\xi-7)/(4\xi-1)}+N^{(8\xi-7)(1-\xi)/(4\xi-1)}N^{2\xi-1} \right) \\
     & \ll
     (MN)^{\varepsilon}N^{(9\xi-6)/(4\xi-1)}.
\end{align*}
   We then deduce that
\begin{align*}
    \mathcal{B}_1(M,
   N)   \ll_{\varepsilon}
   (MN)^{20\varepsilon} \left( M+(MN)^{3/4}+N^{(9\xi-6)/(4\xi-1)} \right).
\end{align*}
  Lemma \ref{lem9} now follows.
\end{proof}

We now proceed to prove Theorem \ref{mainthm}.

\begin{proof}[Proof of Theorem~\ref{mainthm}] Note that it follows from \eqref{initialest} that the estimation given in Lemma \ref{lem9} is valid with
    $\xi=2$.  We further observe that
\begin{align*}
   \frac{3}{2} <  \frac {9\xi-6}{4\xi-1} < \xi,
\end{align*}
  for $\xi > 3/2$ and in the iterative applications of Lemma~\ref{lem9} the exponent of $N$ in the bound for $\mathcal{B}_1(M,N)$ decreases and tends to $3/2$.  We therefore
arrive at the following bound
\begin{align*}
    \mathcal{B}_1(M,
   N)   \ll_{\varepsilon}
   (MN)^{\varepsilon}\left( M+N^{3/2}+(MN)^{3/4} \right),
\end{align*}
   for any $\varepsilon>0$.  Using Lemma \ref{lem2} we then have
\begin{equation*}
\begin{split}
    \mathcal{B}_1(M,
   N)  &  \ll_{\varepsilon}
   (MN)^{\varepsilon}\min \left\{ M+N^{3/2}+(MN)^{3/4},
   N+M^{3/2}+(MN)^{3/4} \right\} \\
   & \ll_{\varepsilon} (MN)^{\varepsilon} \left(M+N+(MN)^{3/4} \right),
   \end{split}
   \end{equation*}
   where the last estimation follows since when $N \leq M, N^{3/2}=N^{3/4}N^{3/4} \leq (MN)^{3/4}$ and similarly when $N \geq M, M^{3/2} \leq (MN)^{3/4}$.
   This establishes Theorem \ref{mainthm}.
   \end{proof}

\section{The Quartic large sieve for Dirichlet Characters}
\label{sec8}

We now proceed to prove Theorem~\ref{quarticlargesieve}.  It is easy to reduce the expression on the left-hand side of
\eqref{final} to a sum of similar expressions with the additional
summation conditions $(q,2)=1$ and $(m,2)=1$ included. Thus it
suffices to estimate
\begin{equation} \label{trans}
\begin{split}
\sum\limits_{\substack{Q<q\le 2Q\\ (q,2)=1}}\
\sideset{}{^\star}\sum\limits_{\substack{\chi \bmod q\\ \chi^4=\chi_0,
\chi^2 \neq \chi_0}} \left| \ \sumstar\limits_{\substack{M<m\le 2M\\
(m,2)=1}} a_m \chi(m)\right|^2 &=
\sumprime\limits_{\substack{n\in \mz[i]\\ Q<\mathcal{N}(n)\le 2Q\\
n\equiv 1 \bmod {(1+i)^3}}} \left| \
\sumstar\limits_{\substack{M< m\le 2M\\ (m,2)=1}} a_m \chi_n(m)\right|^2\\
&= \frac{1}{2}\sumprime\limits_{\substack{n\in \mz[i]\\
Q<\mathcal{N}(n)\le 2Q\\ n\equiv \pm 1 \bmod {(1+i)^3}}} \left| \
\sumstar\limits_{\substack{M<m\le 2M\\ (m,2)=1}} a_m\chi_n(m)
\right|^2,
\end{split}
\end{equation}
where the apostrophe indicates that $n$ is square-free and has no
rational prime divisor and $\chi_n(m)=\left(\frac{m}{n}\right)_4$
is the quartic residue symbol. We shall use this notation for all
$n\in \mz[i]$ and $m\in \mz$.

\subsection{Definition of certain norms}
  In the following, we shall estimate the expression in the last
line of \eqref{trans}. We begin by defining a norm corresponding
to the double sum in the last line of \eqref{trans} by
\begin{equation*}
B_1(Q,M):=\sup\limits_{(a_m)} \| a_m \|^{-2}
\sumprime\limits_{\substack{n \in \mz[i]\\ Q<\mathcal{N}(n)\le 2Q\\
n\equiv \pm 1 \bmod {(1+i)^3}}} \left| \
\sumstar\limits_{\substack{M<m\le 2M\\ (m,2)=1}} a_m
\chi_m(n)\right|^2,
\end{equation*}
where
\begin{equation*}
\|a_m \|^2 = \sum_{m} |a_m|^2,
\end{equation*}
and where by convention we suppose that $(a_m)$ is not identically
zero. \newline

We further define a norm $B_2(Q,M)$ in the same way as $B_1(Q,M)$
except removing the condition that $n$ has no rational prime
divisor. Similarly, we define a norm $B_3(Q,M)$ by further
removing the condition that $n$ is square-free. \newline

  We now use the function $W(x)$ defined in \eqref{W} to see that
\begin{equation*}
  B_3(Q,M) \ll \sup\limits_{(a_m)} \| a_m \|^{-2} \sum\limits_{n\in \mz[i]}
W\left(\frac{\mathcal{N}(n)}{Q}\right) \left| \
\sumstar\limits_{\substack{M<m\le 2M\\ (m,2)=1}} a_m
\chi_m(n)\right|^2.
\end{equation*}
Expanding the sum on the right-hand side, we obtain
\begin{equation} \label{db3bound}
\sumstar\limits_{\substack{M<m_1\le 2M\\ M<m_2\le 2M\\
(m_1m_2,2)=1}} a_{m_1}\overline{a}_{m_2} \sum\limits_{n\in \mz[i]}
W\left(\frac{\mathcal{N}(n)}{Q}\right)
\chi_{m_1}(n)\overline{\chi}_{m_2}(n).
\end{equation}

As in  \cite{DRHB}, it turns out that it suffices to restrict our attention to
the case in which $m_1$ and $m_2$ are coprime.  To see this, we sort the terms
in \eqref{db3bound} according to $\delta = (m_1,m_2)$, and
detecting the condition $(n,\delta)=1$ by the M\"obius function
for $\mz[i]$, we obtain that the expression in \eqref{db3bound}
equals
\begin{equation*}
\begin{split}
\sumstar_{(\delta,2)=1} \sum\limits_{d|\delta} \mu_{[i]}(d)
\sum\limits_{s\in \mz[i]}
W\left(\frac{\mathcal{N}(s)}{Q/\mathcal{N}(d)}\right) & \sumstar\limits_{\substack{M/\delta<r_1\le
2M/\delta \\ M/\delta<r_2\le 2M/\delta\\ (r_1r_2,2\delta)=1\\
(r_1,r_2)=1}} a_{r_1}^*\overline{a_{r_2}^*}
\chi_{r_1}(s)\overline{\chi_{r_2}}(s) \\
& \ll  \sumstar_{(\delta,2)=1} \sum\limits_{d|\delta}
B_4\left(\frac{Q}{\mathcal{N}(d)},\frac{M}{\delta}\right)
\sumstar\limits_{\substack{M/\delta<r\le 2M/\delta \\
(r,2\delta)=1}} |a_{r\delta}|^2,
\end{split}
\end{equation*}
where $d$ runs over non-associate divisors of $\delta$ in
$\mz[i]$,
\[  a_r^*:=a_{r\delta}\chi_r(d), \]
and
\begin{equation*}
B_4(Q,M):=\sup\limits_{(a_m)} \|a_m \|^{-2}
\sumstar\limits_{\substack{M<m_1, m_2\le 2M\\ (m_1m_2,2)=1 \\
(m_1,m_2)=1}} a_{m_1}\overline{a}_{m_2} \sum\limits_{n\in \mz[i]}
W\left(\frac{\mathcal{N}(n)}{Q}\right)
\chi_{m_1}(n)\overline{\chi}_{m_2}(n).
\end{equation*}

Moreover, we define a norm $C_1(M,Q)$ dual to $B_1(Q,M)$ by
\begin{equation*}
C_1(M,Q):=\sup\limits_{(b_{n})} \| b_n \|^{-2}
\sumstar\limits_{\substack{M<m\le 2M\\ (m,2)=1}} \left|
\sumprime\limits_{\substack{n\in \mz[i]\\ Q<\mathcal{N}(n)\le 2Q\\
n\equiv \pm 1 \bmod {(1+i)^3}}} b_{n} \chi_n(m)\right|^2.
\end{equation*}
By the duality principle, we have
\begin{equation} \label{B1C1}
B_1(Q,M)=C_1(M,Q).
\end{equation}
Furthermore, we define a norm $C_2(M,Q)$ by extending the
summation over $m$ in the definition of $C_1(M,Q)$ to all integers
$m$ with $M < m \leq 2M$. Trivially, we have
\begin{equation} \label{C12}
C_1(M,Q)\le C_2(M,Q).
\end{equation}

\subsection{Comparison of the norms} \label{section:comparison}

For the proof of our Theorem \ref{quarticlargesieve}, we need the
following lemma on the norms defined in the previous section.

\begin{lemma} \label{normlemma} Let $Q,M\ge 1$ and $C$ be a sufficiently large positive constant. Then we have the following inequalities:
\begin{equation}
\label{C2e1} C_2(M,Q) \ll (QM)^{\varepsilon}\left(M +
Q^{7/4}\right);
\end{equation}
\begin{equation} \label{C22}
C_2(M,Q)\ll
M^{\varepsilon}Q^{1-1/v}\sum\limits_{j=0}^{v-1}C_2(2^jM^v,Q)^{1/v},
\quad \text{ for each fixed positive integer $v$};
\end{equation}
\begin{equation} \label{B11'}
B_1(Q_1,M)\ll B_1(Q_2,M), \quad \text{ if $Q_1,M\ge 1$ and $Q_2\ge
CQ_1\log (2Q_1M)$};
\end{equation}
\begin{equation} \label{B21}
B_2(Q,M) \ll (\log{2Q})^3 Q^{1/2}X^{-1/2} B_1(XQ^{\varepsilon},M),
\quad \text{for some $X$ with $1\le X\le Q$};
\end{equation}
\begin{equation} \label{B32}
B_3(Q,M) \ll (\log 2Q)^3 Q^{1/2}X^{-1/2} B_2(XQ^{\varepsilon},M),
\quad \text{for some $X$ with $1\le X\le Q$};
\end{equation}
\begin{equation} \label{B34}
B_3(Q,M)\ll M^{\varepsilon}
B_4\left(\frac{Q}{\Delta_1},\frac{M}{\Delta_2}\right), \quad
\text{for some $\Delta_1,\Delta_2\in \mathbb{N}$ with
$\Delta_2^2\ge \Delta_1$};
\end{equation}
\begin{equation} \label{B43}
B_4(Q,M)\ll Q+QM^{\varepsilon-2} \max\left\{B_3(K,M) : K\le
M^4Q^{-1}\right\} + Q^{-1}M^{6+\varepsilon} \sum\limits_{K> M^4/Q}
K^{-2-\varepsilon} B_3(K,M),
\end{equation}
where the sum over $K$ in \eqref{B43} runs over powers of $2$.
\end{lemma}

Since the proofs of \eqref{C22}-\eqref{B43} are essentially the
same as those of (31)-(36) of Lemma 4.1 in \cite{B&Y}, we omit
the proofs here. \newline

We note that it follows from \eqref{B1C1}-\eqref{C22}, that we have
\begin{equation} \label{C2egen}
B_1(Q,M) \ll (QM)^{\varepsilon}\left(Q^{1-1/v} M +
Q^{1+3/(4v)}\right)
\end{equation}
for any $v\in \mn$.

\subsection{Estimating $C_2$} \label{section:C2}

  In this section we prove \eqref{C2e1}. Recall $C_2(M,Q)$
is the norm of the sum
\begin{equation}
\label{eq:C2def} \sum\limits_{M<m\le 2M} \left|
\sumprime\limits_{\substack{n\in \mz[i]\\ Q<\mathcal{N}(n)\le 2Q\\
n\equiv \pm 1 \bmod {(1+i)^3}}}
 b_{n} \chi_n(m)\right|^2,
\end{equation}
where the apostrophe indicates that $n$ is square-free and has no
rational prime divisor. \newline

The sum in \eqref{eq:C2def} is obviously bounded by
\begin{equation}
\label{eq:C2defW}
 \ll \sum\limits_{m\in \mz}
W\left(\frac{m}{M}\right) \left| \sumprime\limits_{n} b_{n}
\chi_n(m)\right|^2,
\end{equation}
where the weight function $W$ is defined as in \eqref{W}.
Expanding out the sum in \eqref{eq:C2defW} we get
\begin{equation*}
\sumprime \limits_{\substack{n_1,n_2\in \mz[i]\\
Q<\mathcal{N}(n_1),\mathcal{N}(n_2)\le 2Q\\ n_1,n_2\equiv \pm 1 \bmod {(1+i)^3}}}
b_{n_1}\overline{b}_{n_2} \sum\limits_{m\in \mz}
W\left(\frac{m}{M}\right) \chi_{n_1}\overline{\chi}_{n_2}(m).
\end{equation*}
Now we extract the greatest common divisor $\Delta$ of $n_1$ and
$n_2$, getting
\begin{equation*}
\sumprime_{\substack{\mathcal{N}(\Delta) \leq 2Q \\ \Delta \equiv 1
\bmod{{(1+i)^3}}}}\ \sumprime \limits_{\substack{n_1,n_2\in
\mz[i]\\ \frac{Q}{\mathcal{N}(\Delta)}<\mathcal{N}(n_1),\mathcal{N}(n_2)\le
\frac{2Q}{\mathcal{N}(\Delta)} \\ n_1,n_2\equiv 1 \bmod {(1+i)^3} \\
(n_1,n_2)=1 \\ (n_1 n_2, \Delta \overline{\Delta}) = 1 }} b_{n_1
\Delta}\overline{b}_{n_2 \Delta} \sum\limits_{\substack{m\in \mz
\\ (m, \mathcal{N}(\Delta)) = 1}} W\left(\frac{m}{M}\right)
\chi_{n_1}\overline{\chi}_{n_2}(m).
\end{equation*}
It is easy to see that there is a one-to-one correspondence between the pairs $(n_1,n_2)$ and $(n'_1\Delta,n'_2 \Delta)$ with $\Delta \equiv 1$, $n_1 \equiv n'_1$ and $n_2 \equiv n'_2 \bmod{(1+i)^3}$.  We write $\delta = (n_1, \bar{n}_2)$ and change variables
via $n_1 \rightarrow \delta n_1$, $n_2 \rightarrow
\overline{\delta} n_2$ to get
\begin{equation} \label{eq:prePoisson}
\begin{split}
\sumprime_{\substack{\mathcal{N}(\Delta) \leq 2Q \\
\Delta \equiv 1 \bmod{{(1+i)^3}}}}\
\sumprime_{\substack{\mathcal{N}(\delta)
\leq \frac{2Q}{\mathcal{N}(\Delta)} \\ \delta \equiv 1 \bmod{{(1+i)^3}} \\
(\mathcal{N}(\delta), \mathcal{N}(\Delta)) = 1}} & \
\sumprime\limits_{\substack{n_1,n_2\in \mz[i]\\
\frac{Q}{\mathcal{N}(\delta \Delta)}<\mathcal{N}(n_1),\mathcal{N}(n_2)\le \frac{2Q}{\mathcal{N}(\delta
\Delta)} \\ n_1,n_2\equiv \pm 1 \bmod {(1+i)^3} \\ (\mathcal{N}(n_1),
\mathcal{N}(n_2)) = 1 \\ (n_1 n_2, \delta \overline{\delta} \Delta
\overline{\Delta}) = 1 } } b_{n_1 \Delta \delta}\overline{b}_{n_2
\Delta \overline{\delta}}
\\
& \times \sum\limits_{\substack{m\in \mz \\ (m, \mathcal{N}(\Delta)) = 1}}
W\left(\frac{m}{M}\right)
\chi_{n_1}\overline{\chi}_{n_2}\overline{\chi}_{\delta}^2(m),
\end{split}
\end{equation}
where we use that for $m \in \mz$,
$\chi_{\delta}\overline{\chi}_{\overline{\delta}}(m)=\chi_{\delta}^2(m)=
\overline{\chi}_{\delta}^2(m)$. Next we remove the coprimality
condition in the sum over $m$ by the M\"obius function, getting
\begin{equation*}
\sum\limits_{\substack{m\in \mz \\ (m, \mathcal{N}(\Delta)) = 1}}
W\left(\frac{m}{M}\right)
\chi_{n_1}\overline{\chi}_{n_2}\overline{\chi}_{\delta}^2(m) =
\sum_{\ell | \mathcal{N}(\Delta)} \mu(\ell)
\chi_{n_1}\overline{\chi}_{n_2}\overline{\chi}_{\delta}^2(\ell)
\sum\limits_{\substack{m\in \mz }} W\left(\frac{m}{M/\ell}\right)
\chi_{n_1}\overline{\chi}_{n_2}\overline{\chi}_{\delta}^2(m),
\end{equation*}
which by the Poisson summation formula is
\begin{equation} \label{afterpoisson}
\begin{split}
\sum_{\ell | \mathcal{N}(\Delta)} \mu(\ell)
\chi_{n_1}\overline{\chi}_{n_2}\overline{\chi}_{\delta}^2(\ell)
& \frac{M}{\ell \mathcal{N}(n_1 n_2 \delta)} \sum\limits_{\substack{h\in \mz}}
\widehat{W}\left(\frac{hM}{\ell \mathcal{N}(n_1 n_2 \delta)}\right)
\\
& \times \sum_{r \bmod{\mathcal{N}(n_1 n_2 \delta)}}
\chi_{n_1}\overline{\chi}_{n_2}\overline{\chi}_{\delta}^2(r)
e\left(\frac{rh}{\mathcal{N}(n_1 n_2 \delta)}\right),
\end{split}
\end{equation}
  where
\begin{align*}
  \widehat{W}(x)=\int\limits^{\infty}_{-\infty}W(y)e(-xy)dy.
\end{align*}

When $h=0$, the expression in \eqref{afterpoisson} vanishes unless
$n_1 = n_2 = \delta =1$. Hence, the contribution of $h=0$ to
\eqref{eq:prePoisson} is
\begin{equation} \label{h0cont}
\ll MQ^{\varepsilon} \sumprime_{\substack{Q<\mathcal{N}(\Delta) \leq 2Q \\
\Delta \equiv 1 \bmod{{(1+i)^3}}}} |b_{\Delta}|^2 \ll M
Q^{\varepsilon}\| b_n \|^2.
\end{equation}

In the sequel, we assume that $h\not=0$. The sum over $r$ in
\eqref{afterpoisson} can be computed by writing $r = r_1 \mathcal{N}(n_2
\delta) + r_2 \mathcal{N}(n_1 \delta) + r_3 \mathcal{N}(n_1 n_2)$ to get
\begin{equation} \label{gauss}
\begin{split}
&\sum_{r \bmod{\mathcal{N}(n_1n_2\delta)}}  \chi_{n_1}\overline{\chi}_{n_2}\overline{\chi}_{\delta}^2(r) e\left(\frac{rh}{\mathcal{N}(n_1 n_2 \delta)}\right)\\
& = \sum_{r_1 \bmod{\mathcal{N}(n_1)}} \chi_{n_1}(r_1\mathcal{N}(n_2\delta))
e\left(\frac{r_1h}{\mathcal{N}(n_1)}\right) \sum_{r_2 \bmod{\mathcal{N}(n_2)}}
\overline{\chi}_{n_2}(r_2\mathcal{N}(n_1\delta))
e\left(\frac{r_2h}{\mathcal{N}(n_2)}\right) \\
&  \hspace*{3cm} \times \sum_{r_3
\bmod{\mathcal{N}(\delta)}}  \overline{\chi}_{\delta}^2(r_3\mathcal{N}(n_1n_2))
e\left(\frac{r_3h}{\mathcal{N}(\delta)}\right) \\
&= \overline{\chi}_{n_1}(h) \chi_{n_2}\chi^2_{\delta}(-h)
\chi_{n_1}(\mathcal{N}(n_2 \delta)) \overline{\chi}_{n_2}(\mathcal{N}(n_1 \delta))
\overline{\chi}_{\delta}^2(\mathcal{N}(n_1n_2)) \tau(\chi_{n_1})
\overline{\tau(\chi_{n_2})\tau(\chi^2_{\delta})},
\end{split}
\end{equation}
   where $\tau (\chi)$ is defined as in \eqref{tau}. \newline

Using quartic reciprocity and the identity
\begin{equation*}
\left(\frac{m}{n}\right)_4 =
\overline{\left(\frac{\overline{m}}{\overline{n}}\right)}_4
\end{equation*}
following from the definition of the quartic residue symbol, we
get the identities
\begin{align*}
\chi_{n}(\mathcal{N}(m))\overline{\chi}_{m}(\mathcal{N}(n))&=\left(\frac{\mathcal{N}(m)}{n}\right)_4
\overline{\left(\frac{\mathcal{N}(n)}{m}\right)}_4 =
\left(\frac{m}{\overline{n}}\right)^2_4 = \chi^2_{\overline{n}}(m)
\\
\chi_{n}(\mathcal{N}(m))\overline{\chi}_{m}^2(\mathcal{N}(n))&=\left(\frac{\mathcal{N}(m)}{n}\right)_4
\left(\frac{\mathcal{N}(n)}{m}\right)^2_4 =
\left(\frac{\mathcal{N}(m)}{\overline{n}}\right)_4 =
\chi_{\overline{n}}(\mathcal{N}(m))
\end{align*}
and
\begin{equation*}
\overline{\chi}_{n}(\mathcal{N}(m))\overline{\chi}_{m}^2(\mathcal{N}(n))=
\overline{\left(\frac{\mathcal{N}(m)}{n}\right)}_4
\overline{\left(\frac{\mathcal{N}(n)}{m}\right)}^2_4 =
\left(\frac{\mathcal{N}(m)}{n}\right)_4 = \chi_{n}(\mathcal{N}(m)),
\end{equation*}
valid for all $m,n\in \mz[i]$ with $m, n \equiv \pm 1
\bmod{{(1+i)^3}}$. We use them to simplify the last line of
\eqref{gauss}, obtaining
\begin{equation*}
  \overline{\chi}_{n_1}(h)
\chi_{n_2}\chi^2_{\delta}(-h)  \left(\frac{n^2_2
\mathcal{N}(\delta)}{\overline{n_1}}\right)_4
\left(\frac{\mathcal{N}(\delta)}{n_2}\right)_4 \tau(\chi_{n_1})
\overline{\tau(\chi_{n_2})\tau(\chi^2_{\delta})}.
\end{equation*}

Now, changing $n_2\rightarrow \overline{n}_2$, the contribution of
$h\not=0$ to the sum in \eqref{afterpoisson} takes the form
\begin{multline}
S_W(M,Q)=M \sumprime_{\substack{\mathcal{N}(\Delta) \leq 2Q \\ \Delta \equiv 1 \bmod{{(1+i)^3}}}}\ \sumprime_{\substack{\mathcal{N}(\delta) \leq
\frac{2Q}{\mathcal{N}(\Delta)} \\ \delta \equiv 1 \bmod{{(1+i)^3}} \\
(\mathcal{N}(\delta), \mathcal{N}(\Delta)) =
1}}\frac{\overline{\tau(\chi_{\delta})}}{\mathcal{N}(\delta)} \sum_{\ell |
\mathcal{N}(\Delta)} \frac{\mu(\ell)}{\ell} \overline{\chi}_{\delta}^2(\ell)
\sum_{h\not=0}  \chi^2_{\delta}(h)
\\
\times \sumprime\limits_{\substack{n_1,n_2\in \mz[i]\\
\frac{Q}{\mathcal{N}(\delta \Delta)}<\mathcal{N}(n_1),\mathcal{N}(n_2)\le \frac{2Q}{\mathcal{N}(\delta
\Delta)} \\ n_1,n_2\equiv \pm 1 \bmod {(1+i)^3} \\ (\mathcal{N}(n_1),
\mathcal{N}(n_2)) = 1 \\ (n_1 n_2, \delta \overline{\delta} \Delta
\overline{\Delta}) = 1 } } \widehat{W}\left(\frac{hM}{\ell \mathcal{N}(n_1 n_2
\delta)}\right) c_{\Delta,\delta,\ell,h,n_1}
c_{\Delta,\delta,\ell,h,n_2}'
\left(\frac{n_1}{n_2}\right)^2_4,\nonumber
\end{multline}
where
\[
c_{\Delta,\delta,\ell,h,n}:=\chi_{n}(\ell)\overline{\chi}_n(h)
\left(\frac{\mathcal{N}(\delta)}{\overline{n}}\right)_4\frac{\tau(\chi_n)}{\mathcal{N}(n)}
b_{n\Delta\delta} \; \; \; \mbox{and} \; \; \;
c_{\Delta,\delta,\ell,h,n}':= \chi_{n}(\ell)\overline{\chi}_n(h)
\left(\frac{\mathcal{N}(\delta)}{\overline{n}}\right)_4\frac{\tau(\chi_n)}{\mathcal{N}(n)}
\overline{b}_{\overline{n}\Delta\overline{\delta}}.
\]
 We now estimate the sum over $n_1$ and $n_2$ directly using
\eqref{eq:quartic}. We denote the inner sum in the definition of
$S_W(M,Q)$ above to be $U(\Delta, \delta, \ell, h)$ so that
\begin{align*}
  U(\Delta, \delta, \ell,
h)= \sumprime\limits_{\substack{n_1,n_2\in \mz[i]\\
\frac{Q}{\mathcal{N}(\delta \Delta)}<\mathcal{N}(n_1),\mathcal{N}(n_2)\le \frac{2Q}{\mathcal{N}(\delta
\Delta)} \\ n_1,n_2\equiv \pm 1 \bmod {(1+i)^3} \\ (\mathcal{N}(n_1),
\mathcal{N}(n_2)) = 1 \\ (n_1 n_2, \delta \overline{\delta} \Delta
\overline{\Delta}) = 1 } } \widehat{W}\left(\frac{hM}{\ell \mathcal{N}(n_1 n_2
\delta)}\right) c_{\Delta,\delta,\ell,h,n_1}
c_{\Delta,\delta,\ell,h,n_2}' \left(\frac{n_1}{n_2}\right)^2_4.
\end{align*}
  To separate the variables $n_1,n_2$, we remove the coprimality
condition $(\mathcal{N}(n_1),\mathcal{N}(n_2))=1$. Because of the
presence of $(\frac{n_1}{n_2})_4$, we may assume that $n_1$ and
$n_2$ are coprime.  Thus $\mathcal{N}(n_1)$ and $\mathcal{N}(n_2)$
must also be coprime unless $n_1$ has a factor in common with
$\overline{n}_2$. We proceed to detect this latter condition using
the Mobius function in the standard way to obtain

\begin{align*}
  U(\Delta, \delta, \ell,
h)&= \sumprime_{\substack{e \in \mz[i] \\ e \equiv 1 \bmod {(1+i)^3} \\ (\mathcal{N}(e), \mathcal{N}(\delta\Delta))=1 }} \mu_{[i]}(e) \sumprime\limits_{\substack{n_1,n_2\in \mz[i]\\
\frac{Q}{\mathcal{N}(\delta
\Delta)}<\mathcal{N}(n_1),\mathcal{N}(n_2)\le
\frac{2Q}{\mathcal{N}(\delta \Delta)} \\ n_1,n_2\equiv \pm 1 \bmod
{(1+i)^3} \\ e|n_1, \overline{e}| n_2 \\ (n_1 n_2, \delta
\overline{\delta} \Delta \overline{\Delta}) = 1 } }
\widehat{W}\left(\frac{hM}{\ell \mathcal{N}(n_1 n_2
\delta)}\right) c_{\Delta,\delta,\ell,h,n_1}
c_{\Delta,\delta,\ell,h,n_2}' \left(\frac{n_1}{n_2}\right)^2_4 \\
&= \sumprime_{\substack{e \in \mz[i] \\ e \equiv 1 \bmod {(1+i)^3}
\\ \mathcal{N}(e) \leq \frac{2Q}{\mathcal{N}(\delta \Delta)} \\
(\mathcal{N}(e), \mathcal{N}(\delta\Delta))=1  }}
\mu_{[i]}(e)\chi_{\mathcal{N}(e)}(\ell)\overline{\chi}_{\mathcal{N}(e)}(h)\left(\frac{\mathcal{N}(\delta)}{\mathcal{N}(e)}\right)_4\frac{\tau(\chi_e)}{\mathcal{N}(e)}
\frac{\tau(\chi_{\overline{e}})}{\mathcal{N}(\overline{e})}\left(\frac{e}{\overline{e}}\right)^2_4 \\
& \hspace*{2cm} \times \sumprime\limits_{\substack{n_1,n_2\in \mz[i]\\
\frac{Q}{\mathcal{N}(e\delta \Delta)}<\mathcal{N}(n_1),\mathcal{N}(n_2)\le \frac{2Q}{\mathcal{N}(e\delta
\Delta)} \\ n_1,n_2\equiv \pm 1 \bmod {(1+i)^3} \\  (n_1 n_2,
e \overline{e}\delta \overline{\delta} \Delta \overline{\Delta}) =
1 } } \widehat{W}\left(\frac{hM}{\ell \mathcal{N}(n_1 n_2 e
\overline{e}\delta)}\right) c_{\Delta,\delta,\ell,h,e,n_1}
c'_{\Delta,\delta,\ell,h,e, n_2} \left(\frac{n_1}{n_2}\right)^2_4 .
\end{align*}
   where
\[ c_{\Delta,\delta,\ell,h,e,n} :=\chi_{n}(\ell)\overline{\chi}_n(h)
\left(\frac{\mathcal{N}(\delta)}{\overline{n}}\right)_4\left(\frac{\mathcal{N}(e)}{n}\right)_4\left(\frac{\mathcal{N}(n)}{e}\right)_4
\left(\frac{n_1}{\overline{e}}\right)^2_4\frac{\tau(\chi_n)}{\mathcal{N}(n)}
b_{ne\Delta\delta} \]
and
\[ c'_{\Delta,\delta,\ell,h,e,n} :=\chi_{n}(\ell)\overline{\chi}_n(h)
\left(\frac{\mathcal{N}(\delta)}{\overline{n}}\right)_4\left(\frac{\mathcal{N}(e)}{n}\right)_4\left(\frac{\mathcal{N}(n)}{\overline{e}}\right)_4\left(\frac{e}{n_2}\right)^2_4\frac{\tau(\chi_n)}{\mathcal{N}(n)}
\overline{b}_{\overline{n}e\Delta\overline{\delta}}. \]
 Next observe that we may freely truncate the sum over $h$ for
\begin{equation*}
|h| \le \frac{Q^2 \ell}{\mathcal{N}(\delta)\mathcal{N}(\Delta)^2M} (QM)^{\varepsilon}=:H
\end{equation*}
since $\widehat{W}$ has rapid decay. More precisely, if we let $S_W
(M,Q) = S'_W(M,Q) + E$ where $S'_W(M,Q)$ is the contribution to
$S_W(M,Q)$ from $0 < |h| \leq H$, then $E \ll
(MQ)^{-100} \| b \|^2$. \newline

  It remains to bound $S'_W(M,Q)$ and we have
\begin{equation*}
\begin{split}
S'_W(M,Q) \ll M \sumprime_{\substack{\mathcal{N}(\Delta) \leq 2Q \\ \Delta
\equiv 1 \bmod{{(1+i)^3}}}}\
\sumprime_{\substack{\mathcal{N}(\delta) \leq
\frac{2Q}{\mathcal{N}(\Delta)} \\ \delta \equiv 1 \bmod{{(1+i)^3}} \\
(\mathcal{N}(\delta), \mathcal{N}(\Delta)) = 1}} & \frac{1}{(\mathcal{N}(\delta))^{1/2}} \sum_{\ell |
\mathcal{N}(\Delta)} \frac{1}{\ell} \\
& \times \sumprime_{\substack{e \in \mz[i] \\ e
\equiv 1 \bmod {(1+i)^3}  \\ \mathcal{N}(e) \leq \frac{2Q}{\mathcal{N}(\delta
\Delta)}\\ (\mathcal{N}(e), \mathcal{N}(\delta\Delta))=1  }}\frac {1}{\mathcal{N}(e)} \sum_{0<
|h| \leq H} \left| U'(\Delta, \delta, \ell, e, h) \right|,
\end{split}
\end{equation*}
  where
\begin{align*}
  U'(\Delta, \delta, \ell, e,
h)=  \sumprime\limits_{\substack{n_1,n_2\in \mz[i]\\
\frac{Q}{\mathcal{N}(e\delta \Delta)}<\mathcal{N}(n_1),\mathcal{N}(n_2)\le \frac{2Q}{\mathcal{N}(e\delta
\Delta)} \\ n_1,n_2\equiv \pm 1 \bmod {(1+i)^3} \\  (n_1 n_2,
e \overline{e}\delta \overline{\delta} \Delta \overline{\Delta}) =
1 } } \widehat{W}\left(\frac{hM}{\ell \mathcal{N}(n_1 n_2 e
\overline{e}\delta)}\right) c_{\Delta,\delta,\ell,h,e,n_1}
c'_{\Delta,\delta,\ell,h,e, n_2} \left(\frac{n_1}{n_2}\right)^2_4 .
\end{align*}
    We  now remove the weight $\widehat{W}$ in $U'(\Delta, \delta, \ell, e,
h)$ by applying Lemma \ref{lem8} to the function $\rho(x) =
\widehat{W}(x)$, which satisfies the conditions of that lemma, as one sees by repeated integration by parts.  We may assume $h>0$ here, since the contribution of the negative $h$'s can be treated similarly and satisfies the same bound.  We use $\sigma = \varepsilon$
to see that
\begin{align*}
    U'(\Delta, \delta, \ell, e,
h) & \ll_{\varepsilon}
\left(\frac{\ell \mathcal{N}(e\overline{e}\Delta)}{hM}\right)^{\varepsilon}\int\limits^{\infty}_{-\infty}|\rho_{+}(\varepsilon+it)||V(\varepsilon+it)|
\dif t,
\end{align*}
  where
\begin{align*}
   V(s)= \sumprime\limits_{\substack{n_1,n_2\in \mz[i]\\
\frac{Q}{\mathcal{N}(e\delta \Delta)}<\mathcal{N}(n_1),\mathcal{N}(n_2)\le \frac{2Q}{\mathcal{N}(e\delta
\Delta)} \\ n_1,n_2\equiv \pm 1 \bmod {(1+i)^3} \\  (n_1 n_2,
e \overline{e}\delta \overline{\delta} \Delta \overline{\Delta}) =
1 } } d_{n_1}d'_{n_2}\left(\frac{n_1}{n_2}\right)^2_4,
\end{align*}
   with
\begin{equation*}
    d_n=c_{\Delta,\delta,\ell,h,e,n}\mathcal{N}(n)^{s} \; \; \;  \mbox{and} \; \; \; d'_n=c'_{\Delta,\delta,\ell,h,e,n}\mathcal{N}(n)^{s} .
\end{equation*}
   Note that $d_{n_1}$ and $d_{n_2}'$ depend on
$\Delta,\delta, \ell, h, e, n$ and $s$, and
\begin{equation*}
|d_n| \ll \left(\frac{\mathcal{N}(e\delta\Delta)}{Q}\right)^{1/2-\varepsilon}
\left| b_{ne\Delta\delta} \right|, \quad |d_n'| \ll
\left(\frac{\mathcal{N}(e\delta\Delta)}{Q}\right)^{1/2-\varepsilon}
\left| b_{\overline{n}e\Delta\overline{\delta}} \right|.
\end{equation*}
 Now, using the Cauchy-Schwarz inequality and the estimate
\eqref{eq:quartic} upon noting that this estimate remains valid if
the summation conditions $m,n\equiv 1 \bmod {(1+i)^3}$ therein are
replaced by $m,n \equiv \pm 1 \bmod {(1+i)^3}$ and
$\left(\frac{n}{m}\right)_4$ replaced by
$\left(\frac{n_1}{n_2}\right)^2_4$, we bound $V(\varepsilon+it)$ by
\begin{equation} \label{square}
\begin{split}
& \left| V(\varepsilon+it) \right|^2 \\
&\ll \sumprime\limits_{\substack{n_1\in \mz[i]\\
\frac{Q}{\mathcal{N}(e\delta \Delta)}<\mathcal{N}(n_1)\le \frac{2Q}{\mathcal{N}(e\delta \Delta)}
\\ n_1 \equiv \pm 1 \bmod {(1+i)^3} \\  (n_1, e
\overline{e}\delta \overline{\delta} \Delta \overline{\Delta}) = 1
} } |d_{n_1}|^2 \;\;\;\; \times
\sumprime\limits_{\substack{n_1\in \mz[i]\\
\frac{Q}{\mathcal{N}(e\delta \Delta)}<\mathcal{N}(n_1)\le \frac{2Q}{\mathcal{N}(e\delta \Delta)}
\\ n_1 \equiv \pm 1 \bmod {(1+i)^3} \\  (n_1, e
\overline{e}\delta \overline{\delta} \Delta \overline{\Delta}) = 1
} } \left|\sumprime\limits_{\substack{n_2\in \mz[i]\\
\frac{Q}{\mathcal{N}(e\delta \Delta)}<\mathcal{N}(n_2)\le \frac{2Q}{\mathcal{N}(e\delta \Delta)}
\\ n_2 \equiv \pm 1 \bmod {(1+i)^3} \\  (n_2, e
\overline{e}\delta \overline{\delta} \Delta \overline{\Delta}) = 1
} } d_{n_2}'
\left(\frac{n_1}{n_2}\right)^2_4 \right|^2 \\
&\ll (QM)^{8\varepsilon}\left(\frac{\mathcal{N}(e\delta
\Delta)}{Q}\right)^{1/4 -4\varepsilon} \left( \sumprime_{\substack{Q/\mathcal{N}(e)< \mathcal{N}(n)\leq 2Q/\mathcal{N}(e)\\
n \equiv \pm 1 \bmod {(1+i)^3}\\ (\mathcal{N}(n), \mathcal{N}(e))=1 }} |b_{en}|^2 \right)^2.
\end{split}
\end{equation}
   Since
\begin{align*}
     \int\limits^{\infty}_{-\infty}|\rho_{+}(\sigma+it)| \dif t \ll_{\varepsilon} 1,
\end{align*}
    we deduce that
\begin{equation} \label{cont}
\begin{split}
S'_W(M,Q) & \ll  M(QM)^{7\varepsilon}
\sumprime_{\substack{\mathcal{N}(\Delta) \leq 2Q \\ \Delta \equiv 1
\bmod{{(1+i)^3}}}}\ \sumprime_{\substack{\mathcal{N}(\delta) \leq
\frac{2Q}{\mathcal{N}(\Delta)} \\ \delta \equiv 1 \bmod{{(1+i)^3}} \\
(\mathcal{N}(\delta), \mathcal{N}(\Delta)) = 1}}\frac{1}{(\mathcal{N}(\delta))^{1/2}} \sum_{\ell |
\mathcal{N}(\Delta)} \frac{1}{\ell} \sumprime_{\substack{e \in \mz[i] \\ e
\equiv 1 \bmod {(1+i)^3}  \\ \mathcal{N}(e) \leq \frac{2Q}{\mathcal{N}(\delta
\Delta)}\\ (\mathcal{N}(e), \mathcal{N}(\delta\Delta))=1  }}\frac {1}{\mathcal{N}(e)}
\\
& \hspace*{2cm} \times \sum_{0< |h| \leq H} \left(\frac{\mathcal{N}(e\delta
\Delta)}{Q}\right)^{\frac14 -2\varepsilon}\sumprime_{\substack{Q/\mathcal{N}(e)< \mathcal{N}(n)\leq 2Q/\mathcal{N}(e)\\
n \equiv \pm 1 \bmod {(1+i)^3}\\ (\mathcal{N}(n), \mathcal{N}(e))=1 }} |b_{en}|^2 \\
& \ll Q^{7/4+2\varepsilon}(QM)^{8\varepsilon}
\sumprime_{\substack{\mathcal{N}(\Delta) \leq 2Q \\ \Delta \equiv 1
\bmod{{(1+i)^3}}}}\frac {1}{\mathcal{N}(\Delta)^{7/4+2\varepsilon}}
\sumprime_{\substack{\mathcal{N}(\delta) \leq
\frac{2Q}{\mathcal{N}(\Delta)} \\ \delta \equiv 1 \bmod{{(1+i)^3}} \\
(\mathcal{N}(\delta), \mathcal{N}(\Delta)) = 1}}\frac{1}{(\mathcal{N}(\delta))^{5/4+2\varepsilon}}
\sum_{\ell | \mathcal{N}(\Delta)}1 \\
& \hspace*{2cm} \times \sumprime_{\substack{e \in \mz[i]
\\ e \equiv 1 \bmod {(1+i)^3}  \\ \mathcal{N}(e) \leq \frac{2Q}{\mathcal{N}(\delta
\Delta)}\\ (\mathcal{N}(e), \mathcal{N}(\delta\Delta))=1  }}\frac
{1}{\mathcal{N}(e)^{3/4+2\varepsilon}}
\sumprime_{\substack{Q/\mathcal{N}(e)< \mathcal{N}(n)\leq 2Q/\mathcal{N}(e)\\
n \equiv \pm 1 \bmod {(1+i)^3}\\ (\mathcal{N}(n), \mathcal{N}(e))=1 }} |b_{en}|^2 \\
& \ll Q^{7/4+2\varepsilon}(QM)^{8\varepsilon}\sumprime_{\substack{Q< \mathcal{N}(n)\leq 2Q\\
n \equiv \pm 1 \bmod {(1+i)^3}}} |b_{n}|^2
\sumprime_{\substack{e \in \mz[i]
\\ e \equiv 1 \bmod {(1+i)^3}  \\  e|n  }}\frac
{1}{\mathcal{N}(e)^{3/4+2\varepsilon}} \\
& \ll Q^{7/4+2\varepsilon}(QM)^{9\varepsilon}\sumprime_{\substack{Q< \mathcal{N}(n)\leq 2Q\\
n \equiv \pm 1 \bmod {(1+i)^3}}} |b_{n}|^2 .
\end{split}
\end{equation}

  Combining \eqref{h0cont} and \eqref{cont}, we deduce that
\eqref{eq:prePoisson} and hence \eqref{eq:C2def} is bounded by
\begin{equation*}
\ll (QM)^{\varepsilon} \left(M+Q^{7/4}\right) \| b_{n} \|^2
\end{equation*}
which implies the desired bound \eqref{C2e1}.

\section{Completion of the proof of Theorem \ref{quarticlargesieve}}
We start with \eqref{C2egen} with any $v\ge 2$ (as one checks
easily that $v=1$ does not lead to any improvement) as an initial
estimate. From \eqref{B21} and \eqref{C2egen}, it follows that
\begin{equation*}
B_2(Q,M) \ll (QM)^{\varepsilon} Q^{1/2}X^{-1/2}
(X^{1+3/(4v)}+X^{1-1/v}M)
\end{equation*}
   for a suitable $X$ with $1\ll X\ll Q$. The worst case is $X = Q$ which shows $B_2(Q,M)$ also satisfies
   \eqref{C2egen}. Repeating the argument, we have
\begin{equation*}
B_3(Q,M)\ll
(QM)^{\varepsilon}\left(Q^{1+3/(4v)}+Q^{1-1/v}M\right).
\end{equation*}
Combining this with \eqref{B43}, we obtain
\begin{eqnarray*}
B_4(Q,M)&\ll& Q+(QM)^{9\varepsilon}QM^{-2} \max\left\{K^{1+3/(4v)}+K^{1-1/v}M \ :\ K\le M^4Q^{-1}\right\}
\\ & & \hspace*{1cm}+ (QM)^{9\varepsilon}M^{6}Q^{-1} \sum\limits_{K\ge M^4/Q} K^{-2-\varepsilon}(K^{1+3/(4v)}+K^{1-1/v}M)\nonumber\\
&\ll&
Q+(QM)^{10\varepsilon}(Q^{-3/(4v)}M^{2+3/v}+Q^{1/v}M^{3-4/v}),
\end{eqnarray*}
where the sum over $K$ runs over powers of $2$. From this and
\eqref{B34}, we deduce that
\begin{equation*}
B_3(Q,M)\ll \frac{Q}{\Delta_1}+(QM)^{\varepsilon}\left(
\left(\frac{Q}{\Delta_1}\right)^{-3/(4v)}
\left(\frac{M}{\Delta_2}\right)^{2+3/v}+
\left(\frac{Q}{\Delta_1}\right)^{1/v}
\left(\frac{M}{\Delta_2}\right)^{3-4/v}\right)
\end{equation*}
for some positive integers $\Delta_1$, $\Delta_2$ with
$\Delta_2^2\ge \Delta_1$. From this and the trivial bound
\begin{equation*}
B_1(Q,M)\ll B_3(Q,M),
\end{equation*}
 we deduce that
\begin{equation}
\label{tem1} B_1(Q,M)\ll
Q+(QM)^{\varepsilon}\left(Q^{-3/(4v)}M^{2+3/v}+Q^{1/v}M^{3-4/v}\right).
\end{equation}
Combining \eqref{tem1} with \eqref{B11'}, we deduce that
\begin{equation} \label{Re1}
B_1(Q,M)\ll
(\tilde{Q}M)^{\varepsilon}\left(\tilde{Q}+\tilde{Q}^{-3/(4v)}M^{2+3/v}+\tilde{Q}^{1/v}M^{3-4/v}\right)
\end{equation}
if $\tilde{Q}\ge CQ\log(2QM)$. We choose $\tilde{Q}:=\max
(Q^{1+\varepsilon}, M^{4-4v/7}). $ Then \eqref{Re1} implies that
\begin{equation} \label{Re2}
B_1(Q,M)\ll (QM)^{\varepsilon}
\left(Q+Q^{1/v}M^{3-4/v}+M^{17/7}\right).
\end{equation}

  It's easy to see that the choice $v=2$ is optimal and a further cycle in the above process does not lead to an improvement of
our result. Combining \eqref{C2egen} with $v=1,2,3$ and
\eqref{Re2} with $v=2$, we obtain our final estimate
\begin{equation}
\label{final0} B_1(Q,M)\ll
(QM)^{\varepsilon}\min\left\{Q^{7/4}+M, \; Q^{11/8}+Q^{1/2}M, \; Q^{5/4}+Q^{2/3}M,\;
Q+Q^{1/2}M+M^{17/7}\right\}.
\end{equation}
which together with \eqref{trans} (noting that the last expression in \eqref{trans} is $\ll B_1(Q,M)$ by the law of quartic reciprocity) implies Theorem \ref{quarticlargesieve}. \hfill $\Box$ \\

Calculating the right-hand side of \eqref{final} for various
ranges of $Q$ and $M$, we obtain that it is bounded by
\begin{equation*}
\ll (QM)^{\varepsilon} \| a_m \|^2
\cdot \left\{ \begin{array}{llll} M &\mbox{ if } Q\le M^{4/7},\\ \\ Q^{7/4} &\mbox{ if } M^{4/7}<Q\le M^{4/5}, \\ \\
Q^{1/2}M &\mbox{ if } M^{4/5}<Q\le M^{8/7}, \\ \\
Q^{11/8} &\mbox{ if } M^{8/7}<Q\le M^{24/17}, \\ \\
Q^{2/3}M &\mbox{ if } M^{24/17}<Q\le M^{12/7}, \\ \\
Q^{5/4} &\mbox{ if } M^{12/7} <Q\le M^{68/35}, \\ \\
M^{17/7} &\mbox{ if } M^{68/35} <Q\le M^{17/7}, \\ \\
Q & \mbox{ if } M^{17/7}<Q.
\end{array} \right.
\end{equation*}
  For convenience, we enclose a table displaying the estimates for $B_1(Q,M)$ that we get
for various ranges. This table should be read as follows. If the
fractions $\alpha$ and $\beta$ are the $(n-1)$-th and $n$-th
entries, respectively, in the first row, and the term $T$ is the
$n$-th entry in the second row, then the estimate $B_1(Q,M) \ll
(QM)^{\varepsilon}T$ holds in the range $M^{\alpha} < Q \leq
M^{\beta}$.

\begin{tabular}{|l||l|l|l|l|l|l|l|l|l|}
\hline Range & 4/7 & 4/5 & 8/7 & 24/17 & 12/7 & 68/35 & 17/7 &  $\infty$ \\
\hline
Bounds & $M$ & $Q^{7/4}$ & $Q^{1/2}$ & $Q^{11/8}$ & $Q^{2/3}$ & $Q^{5/4}$ & $M^{17/7}$ & $Q$  \\
\hline
\end{tabular} \\

It can be easily checked that \eqref{C2egen} with $v\ge 3$ does
not lead to an improvement of \eqref{final0}. \newline

\noindent{\bf Acknowledgments.} During this work, P. G. was
supported by postdoctoral research fellowships at Nanyang
Technological University (NTU) and L. Z. by an AcRF Tier I Grant at NTU.  Moreover, L. Z. was on a visiting stay at {\it Max-Planck-Institut f\"ur Mathematik} (MPIM) in Bonn and he wishes to thank MPIM for its financial support and hospitality during that pleasant and enjoyable period.  Finally, we thank the anonymous referee for his/her many comments and suggestions.

\bibliography{biblio}

\begin{bibdiv}
\begin{biblist}

\bib{Ba2}{article}{
      author={Baier, S.},
       title={The large sieve with quadratic amplitude},
        date={2006},
     journal={Funct. Approx. Comment. Math.},
      volume={36},
       pages={33\ndash 43},
}

\bib{Ba1}{article}{
      author={Baier, S.},
       title={On the large sieve with sparse sets of moduli},
        date={2006},
     journal={J. Ramanujan Math. Soc.},
      volume={21},
       pages={279\ndash 295},
}

\bib{B&Y}{article}{
      author={Baier, S.},
      author={Young, M.~P.},
       title={Mean values with cubic characters},
        date={2010},
     journal={J. Number Theory},
      volume={130},
      number={4},
       pages={879\ndash 903},
}

\bib{SBLZ}{article}{
      author={Baier, S.},
      author={Zhao, L.},
       title={Large sieve inequalities for characters to powerful moduli},
        date={2005},
     journal={Int. J. Number Theory},
      volume={1},
      number={2},
       pages={265\ndash 279},
}

\bib{SBLZ3}{article}{
      author={Baier, S.},
      author={Zhao, L.},
       title={An improvement for the large sieve for square moduli},
        date={2008},
     journal={J. Number Theory},
      volume={128},
      number={1},
       pages={154\ndash 174},
}

\bib{Cla}{article}{
      author={Bombieri, E.},
       title={On the large sieves},
        date={1965},
     journal={Mathematika},
      volume={12},
       pages={201\ndash 225},
}

\bib{BD1}{article}{
      author={Bombieri, E.},
      author={Davenport, H.},
       title={Some inequalities involving trigonometrical polynomials},
        date={1969},
     journal={Ann. Sc. Norm. Super. Pisa},
      volume={23},
       pages={223\ndash 241},
}

\bib{Da}{book}{
      author={Davenport, H.},
       title={{M}ultiplicative {N}umber {T}heory},
     edition={Third edition},
      series={Graduate Texts in Mathematics},
   publisher={Springer-Verlag},
     address={Berlin, etc.},
        date={2000},
      volume={74},
}

\bib{DH1}{article}{
      author={Davenport, H.},
      author={Halberstam, H.},
       title={The values of a trigonometric polynomial at well spaced points},
        date={1966},
     journal={Mathematika},
      volume={13},
       pages={91\ndash 96},
        note={{\it Corrigendum and Addendum}, Mathematika {\bf 14} (1967),
  232-299},
}

\bib{Elliott}{article}{
      author={Elliott, P. D. T.~A.},
       title={On the mean value of $f(p)$},
        date={1970},
     journal={Proc. London Math. Soc. (3)},
      volume={21},
       pages={28\ndash 96},
}

\bib{PXG}{article}{
      author={Gallagher, P.~X.},
       title={The large sieve},
        date={1967},
     journal={Mathematika},
      volume={14},
       pages={14\ndash 20},
}

\bib{DRHB}{article}{
      author={Heath-Brown, D.~R.},
       title={A mean value estimate for real character sums},
        date={1995},
     journal={Acta Arith.},
      volume={72},
      number={3},
       pages={235\ndash 275},
}

\bib{DRHB1}{article}{
      author={Heath-Brown, D.~R.},
       title={Kummer's conjecture for cubic {G}auss sums},
        date={2000},
     journal={Israel J. Math.},
      volume={120},
       pages={97\ndash 124},
}

\bib{H&P}{article}{
      author={Heath-Brown, D.~R.},
      author={Patterson, S.~J.},
       title={The distribution of {K}ummer sums at prime arguments},
        date={1979},
     journal={J. Reine Angew. Math.},
      volume={310},
       pages={111\ndash 130},
}

\bib{I&R}{book}{
      author={Ireland, K.},
      author={Rosen, M.},
       title={{A} {C}lassical {I}ntroduction to {M}odern {N}umber {T}heory},
     edition={Second edition},
      series={Graduate Texts in Mathematics},
   publisher={Springer-Verlag},
     address={New York},
        date={1990},
      number={84},
}

\bib{JVL1}{article}{
      author={Linnik, J.~V.},
       title={The large sieve},
        date={1941},
     journal={Doklady Akad. Nauk SSSR},
      volume={36},
       pages={119\ndash 120},
}

\bib{Luo}{article}{
      author={Luo, W.},
       title={On {H}ecke {$L$}-series associated with cubic characters},
        date={2004},
     journal={Compos. Math.},
      volume={140},
      number={5},
       pages={1191\ndash 1196},
}

\bib{Lem}{book}{
      author={Montgomery, H.~L.},
       title={{T}opics in {M}ultiplicative {N}umber {T}heory},
      series={Lecture Notes in Mathematics},
   publisher={Spring-Verlag},
     address={Barcelona, etc.},
        date={1971},
      volume={227},
}

\bib{HM}{book}{
      author={Montgomery, H.~L.},
       title={{T}opics in {M}ultiplicative {N}umber {T}heory},
      series={Lecture Notes in Mathematics},
   publisher={Spring-Verlag},
     address={Berlin, etc.},
        date={1971},
      volume={227},
}

\bib{HM2}{article}{
      author={Montgomery, H.~L.},
       title={The analytic priciples of large sieve},
        date={1978},
     journal={Bull. Amer. Math. Soc.},
      volume={84},
      number={4},
       pages={547\ndash 567},
}

\bib{MVa}{article}{
      author={Montgomery, H.~L.},
      author={Vaughan, R.~C.},
       title={The large sieve},
        date={1973},
     journal={Mathematika},
      volume={20},
       pages={119\ndash 134},
}

\bib{Patt}{article}{
      author={Patterson, S.~J},
       title={The distribution of general {G}auss sums and similar arithmetic
  functions at prime arguments},
        date={1987},
     journal={Prod. London Math. Soc. (3)},
      volume={54},
      number={2},
       pages={193\ndash 215},
}

\bib{PrRa}{article}{
      author={Prakash, G.},
      author={Ramana, D.~S.},
       title={The large sieve with quadratic amplitude},
        date={2009},
     journal={J. Number Theory},
      volume={129},
      number={2},
       pages={428\ndash 433},
}

\bib{Wol}{article}{
      author={Wolke, D.},
       title={On the large sieve with primes},
        date={1971/72},
     journal={Acta Math. Acad. Sci. Hungar.},
      volume={22},
       pages={239\ndash 247},
}

\bib{Zha4}{article}{
      author={Zhao, L.},
       title={Large sieve inequalities for special characters to prime square
  moduli},
        date={2004},
     journal={Funct. Approx. Comment. Math.},
      volume={32},
       pages={99\ndash 106},
}

\bib{Zha}{article}{
      author={Zhao, L.},
       title={Large sieve inequality for characters to square moduli},
        date={2004},
     journal={Acta Arith.},
      volume={112},
      number={3},
       pages={297\ndash 308},
}

\bib{Zha3}{article}{
      author={Zhao, L.},
       title={An improvement of a large sieve inequality in high dimensions},
        date={2005},
     journal={Mathematika},
      volume={52},
       pages={93\ndash 100},
}

\bib{Zha2}{article}{
      author={Zhao, L.},
       title={Large sieve inequalities with quadratic amplitudes},
        date={2007},
     journal={Monatsh. Math.},
      volume={151},
      number={2},
       pages={165\ndash 173},
}

\end{biblist}
\end{bibdiv}
\bibliographystyle{amsxport}

\vspace*{.5cm}

\noindent\begin{tabular}{p{8cm}p{8cm}}
Peng Gao & Liangyi Zhao \\
Div. of Math. Sci., School of Phys. \& Math. Sci., & Div. of Math. Sci., School of Phys. \& Math. Sci., \\
Nanyang Technological Univ., Singapore 637371 & Nanyang Technological Univ., Singapore 637371 \\
Email: {\tt penggao@ntu.edu.sg} & Email: {\tt lzhao@pmail.ntu.edu.sg} \\
\end{tabular}

\end{document}